# Stabilization and Spill-Free Transfer of Viscous Liquid in a Tank

**Iasson Karafyllis*** **and Miroslav Krstic****

*Dept. of Mathematics, National Technical University of Athens,
Zografou Campus, 15780, Athens, Greece
email: iasonkar@central.ntua.gr

**Dept. of Mechanical and Aerospace Eng., University of California, San Diego,
La Jolla, CA 92093-0411, U.S.A.
email: krstic@ucsd.edu

## Abstract

Flow control occupies a special place in the fields of partial differential equations (PDEs) and control theory, where the complex behavior of solutions of nonlinear dynamics in very high dimension is not just to be understood but also to be assigned specific desired properties, by feedback control. Among several benchmark problems in flow control, the liquid-tank problem is particularly attractive as a research topic. It is among the few hard problems in PDE control where the solutions are readily comprehensible—with feedback laws that are relatively simple in form, with a clear physical meaning.

In the liquid-tank problem the objective is to move a tank filled with liquid, suppress the nonlinear oscillations of the liquid in the process, bring the tank and liquid to rest, and avoid liquid spillage in the process. In other words, this is a problem of nonlinear PDE stabilization subject to state constraints.

This review article focuses only on recent results on liquid-tank stabilization for viscous (incompressible, Newtonian) liquids, as viscosity makes the earlier control designs for inviscid Saint-Venant models inapplicable. All possible cases are studied: with and without friction from the tank walls, with and without surface tension. Moreover, results that have not been previously published are provided for the linearization of the tank-liquid system. The linearization of the tank-liquid system gives a high-order PDE which is a combination of a wave equation with Kelvin-Voigt damping and an Euler-Bernoulli beam equation.

The feedback design methodology presented in the article is based on Control Lyapunov Functionals (CLFs), suitably extended from the CLF methodology for ODEs to the infinite-dimensional case. The CLFs proposed are modifications and augmentations of the total energy functionals for the tank-liquid system, so that the dissipative effects of viscosity, friction, and surface tension are captured and additional dissipation by feedback is made relatively easy.

The article closes with an extensive list of open problems.

## Contents

## 1. Introduction

***Broader Context: Flow Control.*** For a decade or more in the late 1990s and early 2000s, the Navier-Stokes and other partial differential equation (PDE) models of fluid flows were among the highly attractive research topics in applied mathematics and theoretical engineering. Model reduction, finite-dimensional attractors, hydrodynamic stability studies of transition to turbulence, and flow control were some of the topics that drew hundreds of researchers and through which fluid dynamics, as a physics discipline, intersected with mathematics and control theory.

Among these mathematical topics in fluid dynamics particularly notable was flow control, which aimed not only to understand the solutions to fluid flow PDEs but to also influence them. With approaches ranging from controllability (open-loop controls) to optimal control (open- and closed-loop controls) to feedback design for Lyapunov stabilization (closed-loop controls), a considerable body of knowledge was built, captured in several books, such as, for example, [1, 13, 24, 54], and hundreds of papers. It would take a very long survey to review all these results, on finite or infinite time horizons, local or global, valid for small or unlimited Reynolds numbers, for laminar or turbulent regimes, using actuation on the boundary or within the domain, with normal or tangential velocity inputs at the boundary, in 2-D or 3-D, and in a variety of spatial geometries, some simple or periodic, and some complex and thus seemingly realistic but not universally representative of the dynamics of the fluid, especially in the turbulent regime.

***Tank-Liquid System: Application, Modeling, and Control.*** In this review we focus on one particular canonical topic in flow control—the stabilization of a liquid in a tank being transferred to a desired position of the tank. To a mathematician and a control theorist this problem is interesting for two reasons. First, the system is nonlinear and has two spatially-distributed infinite-dimensional states (liquid height and velocity). Second, this (doubly) infinite-dimensional system has only one input—the force acting on the tank. Hence, this physically relevant problem also represents a benchmark in nonlinear PDE control in the "least actuated" (hardest) sense—with merely one input.

The liquid-tank system is not only of theoretical interest. The problem has an intensely practical motivation. For instance, in [57, 58] on learns of challenges that arise in spacecraft operation due to the liquid fuel sloshing during spacecraft transfer and maneuvering, due to the dynamic interaction between the liquid and the rigid body. One way of mitigating this interaction is through a design of devices that can achieve sloshing suppression (see for instance Chapter 3 in [28]). However, suppression of the unsteady motion of the liquid by active control is superior in achieving suppression because it requires no additional devices. The actuators that perform the transfer, and thus induce the fluid-body unsteady interaction, are available also for suppressing this interaction.

From a mathematical point of view, the description of sloshing is a highly non-trivial modeling problem. Two main approaches have been used for the modeling of free surface flows of incompressible liquids in the literature: (i) the use of the fluid momentum equations under the assumption of the irrotational flow of the liquid (see for instance [28, 29, 41]), and (ii) the use of the fluid momentum equations for the derivation of Saint-Venant models (see [2]-the first paper by Adhémar Jean Claude Barré de Saint-Venant in 1871) by neglecting the fluid motion in the direction of the liquid height. The Saint-Venant model or shallow water model is a well-known mathematical model which has been used extensively and many modifications of this model take into account various types of forces like gravity, viscous stresses, surface tension, friction forces (see [6, 7, 8, 11, 22, 39, 40, 41, 44, 45, 53]). In this review we focus on sloshing induced by the movement of the container and on 1-D Saint-Venant models for the description of the liquid motion.

Control studies of the Saint-Venant model have focused on the inviscid Saint-Venant model (i.e., the model that ignores viscous stresses and surface tension) and its linearization around an equilibrium point (see [3, 4, 5, 12, 13, 14, 15, 17, 18, 19, 20, 21, 25, 26, 42, 47, 48]). The feedback design has been performed by employing either the backstepping methodology or the Control Lyapunov Functional (CLF) methodology. Although there is no inviscid liquid, the use of the inviscid Saint-Venant model is justified when studying the flow in rivers: in this case the inertial, gravity and friction forces are orders of magnitude larger than the viscous stresses and the effect of viscosity is negligible.

However, when one studies the flow in a tank there is no guarantee that the effect of viscosity is negligible because the velocity of the fluid is (expected to be) relatively small. To this purpose, viscous Saint-Venant models have been proposed and studied in [6, 7, 22, 33, 34, 35, 36, 39, 45, 53]. From a mathematical point of view the effect of the viscosity is huge: in the case where surface tension is absent the system is described by two Ordinary Differential Equations (ODEs), one first-order hyperbolic Partial Differential Equation (PDE) and one parabolic PDE, whereas in the inviscid case we have two ODEs and two first-order hyperbolic PDEs.

From the point of view of applications, when one studies the flow in a tank, the avoidance of the phenomenon of liquid spilling out of the tank becomes as important as the sloshing problem. Thus, the solution of the spill-free and slosh-free movement problem by means of a robust feedback law becomes a significant mathematical problem with possibly important applications. The recent works [33, 34, 35, 36] study this particular feedback stabilization problem for the viscous Saint-Venant model without linearization around an equilibrium point. More specifically, in [33, 34] a feedback control law is constructed by employing the Control Lyapunov Functional (CLF) methodology for the viscous Saint-Venant model without wall friction and surface tension (state feedback in [33] and output feedback in [34]). In [35, 36] it is shown that the same feedback control law proposed in [33] works even if friction forces or surface tension are present. It should be noticed that [33, 34, 35, 36] are the only works that guarantee a spill-free movement of the fluid.

The key element of the control design methodology that we review is the construction of the CLF, which then make the choice of the stabilizing feedback law relatively immediate, though certainly not trivial. The construction of the CLFs is among the most interesting ideas we introduce in this article. The CLFs are built on physical concepts of energy for the tank-liquid system, but the CLFs are never simply the energy functionals. They contain modifications and augmentations of the energy functionals, which capture the energy (and possibly dissipating) effects of physical phenomena like viscosity, friction, gravity or surface tension, and allow the dissipation to be enhanced by feedback based on the CLF.

Let us also point out the obvious—the tank-liquid problem is not just a problem in flow control but actually an elementary form of a problem in fluid-structure interaction, where the structure simply the rigid tank. One can, therefore, regard the problem of control of the tank-liquid system as a first step towards the development of controllers to attenuate fluid-structure interaction, such as flutter. As the reader will note, the controllers we design employ feedback from both the state of the liquid (height and velocity profiles) and the state of the tank (position and velocity).

***Linearized Tank-Liquid System.*** After we review, in the article, the existing results for the spill-free and slosh-free movement problem of the 1-D tank-liquid system, we proceed to study some of the same problems for the linearized version of the tank-liquid system. The linearization gives a high-order PDE which is a combination of a wave equation with Kelvin-Voigt damping and an Euler-Bernoulli beam equation. There are important differences between the results for the linearization and the results for the nonlinear tank-liquid system:

1) In the linearized case we provide existence/uniqueness results for the closed-loop system. In the nonlinear case we do not provide existence/uniqueness results for the corresponding closed-loop system.

2) In the linearized case we can study the situation where both friction and surface tension are present. In the nonlinear case, we cannot study the situation where both friction and surface tension are present.

3) The state norm for which stabilization is achieved in the linearized case is stronger than the state norm for which stabilization is achieved in the nonlinear case. This difference is explained by the difference of the CLFs in the nonlinear and the linearized case.

The paper is structured as follows. In Section 2 we describe the mathematical ODE-PDE model of the tank-liquid system. In Section 3 we present the control problem and we review the existing results for the spill-free and slosh-free movement problem of the 1-D tank-liquid system. In Section 4 we give the linearization of the system as well as some important properties of the open-loop system. Section 5 provides the stabilization results for the linearized model and a detailed comparison with the results for the nonlinear case. Section 6 is devoted to the proofs of all results in the paper. Finally, in Section 7 we give an extensive list of important problems that remain open and can be topics for future research.

***Notation.*** Throughout the article we adopt the following notation.

* $\mathbb{R}_+ = [0,+\infty)$ denotes the set of non-negative real numbers.

* Let $S \subseteq \mathbb{R}^n$ be an open set and let $A \subseteq \mathbb{R}^n$ be a set that satisfies $S \subseteq A \subseteq cl(S)$. By $C^0(A;\Omega)$, we denote the class of continuous functions on $A$, which take values in $\Omega \subseteq \mathbb{R}^m$. By $C^k(A;\Omega)$, where $k \ge 1$ is an integer, we denote the class of functions on $A \subseteq \mathbb{R}^n$, which takes values in $\Omega \subseteq \mathbb{R}^m$ and has continuous derivatives of order $k$. In other words, the functions of class $C^k(A;\Omega)$ are the functions which have continuous derivatives of order $k$ in $S = \mathrm{int}(A)$ that can be continued continuously to all points in $\partial S \cap A$. When $\Omega = \mathbb{R}$ then we write $C^0(A)$ or $C^k(A)$. When $I \subseteq \mathbb{R}$ is an interval and $G \in C^1(I)$ is a function of a single variable, $G'(h)$ denotes the derivative with respect to $h \in I$.

* Let $I \subseteq \mathbb{R}$ be an interval and let $Y$ be a normed linear space. By $C^0(I;Y)$, we denote the class of continuous functions on $I$, which take values in $Y$. By $C^1(I;Y)$, we denote the class of continuously differentiable functions on $I$, which take values in $Y$.

* Let $I \subseteq \mathbb{R}$ be an interval, let $a < b$ be given constants and let $u: I \times [a,b] \to \mathbb{R}$ be a given function. We use the notation $u[t]$ to denote the profile at certain $t \in I$, i.e., $(u[t])(x) = u(t,x)$ for all $x \in [a,b]$. When $u(t,x)$ is (twice) differentiable with respect to $x \in [a,b]$, we use the notation $u_x(t,x)$ ($u_{xx}(t,x)$) for the (second) derivative of $u$ with respect to $x \in [a,b]$, i.e., $u_x(t,x) = \frac{\partial u}{\partial x}(t,x)$ ($u_{xx}(t,x) = \frac{\partial^2 u}{\partial x^2}(t,x)$). When $u(t,x)$ is differentiable with respect to $t$, we use the notation $u_t(t,x)$ for the derivative of $u$ with respect to $t$, i.e., $u_t(t,x) = \frac{\partial u}{\partial t}(t,x)$. When $u[t] \in X$ for all $t \in I$, where $X$ is a normed linear space with norm $\| \ \|_X$ and the mapping $I \ni t \to u[t] \in X$ is $C^1$, i.e., the exists a continuous mapping $v: I \to X$ with $\lim_{h \to 0}\left( \left\| \frac{u[t+h]-u[t]}{h} - v[t] \right\|_X \right) = 0$ for all $t \in I$, we use the notation $u_t$ for the mapping $v: I \to X$. Furthermore, when $u \in C^1(I;X)$ and the mapping $I \ni t \to u_t[t] \in X$ is $C^1$, i.e., the exists a continuous mapping $w: I \to X$ with $\lim_{h \to 0}\left( \left\| \frac{u_t[t+h]-u_t[t]}{h} - w[t] \right\|_X \right) = 0$ for all $t \in I$, we use the notation $u_{tt}$ for the mapping $w: I \to X$. Mixed derivatives are to be understood in this way. For example, when $u_x \in C^1(I;X)$, i.e., the exists a continuous mapping $\varphi: I \to X$ with $\lim_{h \to 0}\left( \left\| \frac{u_x[t+h]-u_x[t]}{h} - \varphi[t] \right\|_X \right) = 0$ for all $t \in I$, we use the notation $u_{xt}$ for the mapping $\varphi: I \to X$.

* Given a set $U \subseteq \mathbb{R}^n$, $\chi_U$ denotes the characteristic function of $U$, i.e. the function defined by $\chi_U(x) := 1$ for all $x \in U$ and $\chi_U(x) := 0$ for all $x \notin U$.

* Let $a<b$ be given constants. For $p\in[1,+\infty)$, $L^p(a,b)$ is the set of equivalence classes of Lebesgue measurable functions $u:(a,b)\to\mathbb{R}$ with $\|u\|_p := \left(\int_a^b |u(x)|^p\,dx\right)^{1/p} < +\infty$. The scalar product in $L^2(a,b)$ is denoted by $\langle \bullet,\bullet \rangle$, i.e., $\langle f,g\rangle = \int_a^b f(x)g(x)dx$ for all $f,g\in L^2(a,b)$. $L^\infty(a,b)$ is the set of equivalence classes of Lebesgue measurable functions $u:(a,b)\to\mathbb{R}$ with $\|u\|_\infty := \underset{x\in(a,b)}{ess\sup}\left(|u(x)|\right) < +\infty$. For an integer $k\geq 1$, $H^k(a,b)$ denotes the Sobolev space of functions in $L^2(a,b)$ with all its weak derivatives up to order $k\geq 1$ in $L^2(a,b)$.

## 2. The Mathematical Model

We consider a one-dimensional model for the motion of a tank. The tank contains a viscous, Newtonian, incompressible liquid. The tank is subject to a force that can be manipulated. We assume that the liquid pressure is hydrostatic and consequently, the liquid is modeled by the one-dimensional (1-D) viscous Saint-Venant equations, whereas the tank obeys Newton's second law and consequently we consider the tank acceleration to be the control input.

### *2.1. A General 1-D Model*

We next give a general 1-D model for the tank-liquid system that takes into account all possible forces exerted on the fluid: gravity, viscous stresses, surface tension and friction. Let the position of the left side of the tank at time $t\geq 0$ be $a(t)$ and let the length of the tank be $L>0$ (a constant). The equations describing the motion of the liquid within the tank are

$$H_t + (H\bar{v})_z = 0\text{, for } t>0,\ z\in\left[a(t),a(t)+L\right] \tag{2.1}$$

$$\begin{aligned}(H\bar{v})_t + \left(H\bar{v}^2 + \frac{1}{2}gH^2\right)_z - \sigma H\left(\frac{H_{zz}}{\left(1+H_z^2\right)^{3/2}}\right)_z &, \\ = \mu\left(H\bar{v}_z\right)_z - \kappa\left(H(t,z),\bar{v}(t,z)-\dot{a}(t)\right)\left(\bar{v}(t,z)-\dot{a}(t)\right) & \\ \text{for } t>0,\ z\in\left(a(t),a(t)+L\right) &\end{aligned} \tag{2.2}$$

where $H(t,z)>0$, $\bar{v}(t,z)\in\mathbb{R}$ are the liquid level and the liquid velocity, respectively, at time $t\geq 0$ and position $z\in\left[a(t),a(t)+L\right]$, $\kappa\in C^0\left((0,+\infty)\times\mathbb{R};\mathbb{R}_+\right)$ is the friction coefficient that depends on the liquid level and the relative velocity of the fluid with respect to the tank, while $g,\mu>0$, $\sigma\geq 0$ (constants) are the acceleration of gravity, the kinematic viscosity of the liquid and the ratio of the surface tension and liquid density, respectively. In certain works the term $\left(\frac{H_{zz}}{\left(1+H_z^2\right)^{3/2}}\right)_z$ is replaced by $H_{zzz}$ (see [6, 7, 8, 44]), but here we use a more accurate description of the surface tension. Equations (2.1), (2.2) can be derived by performing mass and momentum balances (from first principles assuming that the liquid pressure is the combination of hydrostatic pressure and capillary pressure given by the Young-Laplace equation-see [16]).

Various empirical relations have been used for the friction coefficient in the literature:

- in [5, 25] the authors use the relation $\kappa(h,v) = c_f |v|$, where $c_f > 0$ is a constant,
- in [7] the authors use the relation $\kappa(h,v) = r_0 + r_1 h|v|$, where $r_0 > 0$, $r_1 \geq 0$ are constants,
- in [20] the authors use the relation $\kappa(h,v) = r_2 h^{-1/3} \left(b_2 + 2h\right)^{4/3} |v|$, where $r_2 > 0$ is a constant and $b_2 > 0$ is the (constant) width of the tank,
- in [22] the authors derive the velocity-independent relation $\kappa(h) = 3\mu b_3 / \left(3\mu + 4b_3 h\right)$, where $b_3 > 0$ is a constant.

The liquid velocities at the walls of the tank must coincide with the tank velocity, i.e., we have:

$$\bar{v}(t,a(t)) = \bar{v}(t,a(t)+L) = w(t) \text{, for } t \geq 0 \tag{2.3}$$

where $w(t) = \dot{a}(t)$ is the velocity of the tank at time $t \geq 0$. Moreover, since the tank acceleration is the control input, we get

$$\ddot{a}(t) = -f(t) \text{, for } t > 0 \tag{2.4}$$

where $-f(t)$, the control input to the problem, is equal to the force exerted on the tank at time $t \geq 0$ divided by the total mass of the tank. Using (2.1) and (2.3), it becomes clear that every classical solution of (2.1), (2.3) satisfies $\frac{d}{dt}\left( \int_{a(t)}^{a(t)+L} H(t,z)dz \right) = 0$ for all $t > 0$. Therefore, the total mass of the liquid is constant. Thus, without loss of generality, we assume that the following equation holds

$$\int_{a(t)}^{a(t)+L} H(t,z)dz \equiv m \tag{2.5}$$

where $m > 0$ is the total mass of the liquid divided by the product of liquid density times the width of the tank. It should be emphasized that for obvious physical reasons, the liquid level $H(t,z)$ must be positive for all times, i.e., we must have:

$$\min_{x \in [0,L]} \left( H(t,a(t)+x) \right) > 0 \text{, for } t \geq 0 \tag{2.6}$$

For a complete mathematical model of the system in the case $\sigma > 0$ (the case where surface tension is present), we need two additional boundary conditions that describe the interaction between the liquid and the solid walls of the tank. There are many ways to describe the evolution of the angle of contact of a liquid with a solid boundary (see the detailed presentation in [38]). In [49, 50], Schweizer suggested (based on energy arguments and the fact that there may be a discrepancy between the actual microscopic and the apparent macroscopic contact angle) the use of a constant contact angle. Moreover, the assumption of a constant contact angle allows the well-posedness of the overall problem (at least for small data; see [49, 50, 55]). The constant contact angle approach has been used extensively in the literature (see for instance [27, 55, 56]). In this work, we adopt the constant contact angle approach by imposing a contact angle equal to $\pi / 2$. Therefore, the model is accompanied by the following boundary conditions (written in a way that holds even in the case $\sigma = 0$, i.e., the case where surface tension is absent and the additional boundary conditions are not needed):

$$\sigma H_z(t,a(t)) = \sigma H_z(t,a(t)+L) = 0 \text{, for } t \ge 0 \tag{2.7}$$

Applying the transformation

$$\begin{aligned} v(t,x) &= \bar{v}(t,a(t)+x) - w(t) \\ h(t,x) &= H(t,a(t)+x) \\ \xi(t) &= a(t) - a^* \end{aligned} \tag{2.8}$$

where $a^* \in \mathbb{R}$ is the specified position (a constant) to which we want to bring (and maintain) the left side of the tank, we obtain the model:

$$\dot{\xi} = w \quad , \quad \dot{w} = -f \text{, for } t \ge 0 \tag{2.9}$$

$$h_t + (hv)_x = 0 \text{, for } t > 0,\ x \in [0, L] \tag{2.10}$$

$$(hv)_t + \left( hv^2 + \frac{1}{2} g h^2 \right)_x - \sigma h \left( \frac{h_{xx}}{\left(1+h_x^2\right)^{3/2}} \right)_x = \mu \left( h v_x \right)_x - \kappa(h,v)v + hf \text{, for } t > 0,\ x \in (0, L) \tag{2.11}$$

$$v(t,0) = v(t,L) = 0 \text{, for } t \ge 0 \tag{2.12}$$

$$\int_0^L h(t,x)dx = m \text{, for } t \ge 0 \tag{2.13}$$

$$\sigma h_x(t,0) = \sigma h_x(t,L) = 0 \text{, for } t \ge 0 \tag{2.14}$$

$$\min_{x \in [0,L]} \left( h(t,x) \right) > 0 \text{, for } t \ge 0 \tag{2.15}$$

The open-loop system (2.9)-(2.15), i.e., system (2.9)-(2.15) with $f(t) \equiv 0$, allows a continuum of equilibrium points, namely the points

$$h(x) \equiv h^*,\ v(x) \equiv 0 \text{, for } x \in [0, L] \tag{2.16}$$

$$\xi \in \mathbb{R},\ w = 0 \tag{2.17}$$

where $h^* = m / L$. The existence of a continuum family of equilibrium points for the open-loop system given by (2.16), (2.17), with the family parametrized by an arbitrary position of the tank, while the liquid is at a unique and spatially constant height, implies that the desired equilibrium point, i.e., the equilibrium point with $\xi = 0$, is not asymptotically stable for the open-loop system.

## *2.2. The Liquid-Tank Interaction*

Equations (2.2), (2.3) and (2.7) allow us to calculate the force $F_{\text{liquid}}(t)$ exerted on the liquid divided by the product of liquid density times the width of the tank:

$$\begin{aligned} F_{\text{liquid}}(t) &= \mu \left( H(t,a(t)+L)\bar{v}_z(t,a(t)+L) - H(t,a(t))\bar{v}_z(t,a(t)) \right) \\ &- \int_{a(t)}^{a(t)+L} \kappa\left( H(t,z), \bar{v}(t,z) - \dot{a}(t) \right) \left( \bar{v}(t,z) - \dot{a}(t) \right) dz \\ &+ \frac{g}{2} \left( H^2(t,a(t)) - H^2(t,a(t)+L) \right) \\ &+ \sigma \left( H(t,a(t)+L) H_{zz}(t,a(t)+L) - H(t,a(t)) H_{zz}(t,a(t)) \right) \end{aligned} \tag{2.18}$$

Assuming that the tank is subject to an external force $F_{\text{external}}(t)$ (divided by the product of liquid density times the width of the tank) and using the transformation (2.8), equation (2.18), Newton's laws and (2.4) we obtain the equation:

$$
\begin{aligned}
&\bar{m}f(t)=\mu\big(h(t,L)v_x(t,L)-h(t,0)v_x(t,0)\big)\\
&-\int_0^L \kappa\big(h(t,x),v(t,x)\big)v(t,x)dx+\frac{g}{2}\big(h^2(t,0)-h^2(t,L)\big)\\
&+\sigma\big(h(t,L)h_{xx}(t,L)-h(t,0)h_{xx}(t,0)\big)-F_{\text{external}}(t)
\end{aligned}
\tag{2.19}
$$

where $\bar{m}>0$ is the tank mass divided by the product of liquid density times the width of the tank.

Equation (2.19) implies that the absence of an external force on the tank does not necessarily imply that the tank acceleration $-f(t)$ is zero. Indeed, the tank may move due to the tank-liquid interaction. Equation (2.19) is also important when the actual control input is the external force on the tank $F_{\text{external}}(t)$. Indeed, using (2.10), (2.12) and (2.19) we can relate directly the external force on the tank and the tank acceleration:

$$
\begin{aligned}
&F_{\text{external}}(t)=-\int_0^L \kappa\big(h(t,x),v(t,x)\big)v(t,x)dx\\
&+\mu\frac{d}{dt}\big(h(t,0)-h(t,L)\big)+\frac{g}{2}\big(h^2(t,0)-h^2(t,L)\big)\\
&+\sigma\big(h(t,L)h_{xx}(t,L)-h(t,0)h_{xx}(t,0)\big)-\bar{m}f(t)
\end{aligned}
\tag{2.20}
$$

Equation (2.20) may be used in a straightforward way for the calculation of the external force that must be applied on the tank in order to have a tank acceleration equal to $-f(t)$.

## 3. A Review of Stabilization Results

The control objective is to drive asymptotically the tank to a specified position without liquid spilling out and having both the tank and the liquid within the tank at rest. From a mathematical point of view our objective is to design a (potentially) full-state feedback law of the form

$$
f(t)=\bar{F}\big(\xi(t),w(t),h[t],v[t]\big),\text{ for } t\geq 0, \tag{3.1}
$$

which achieves stabilization of the equilibrium point with $\xi=0$. Moreover, we additionally require that the following "no-spill condition" holds

$$
\max_{x\in[0,L]}\big(h(t,x)\big)<H_{\max},\text{ for } t\geq 0 \tag{3.2}
$$

where $H_{\max}>0$ is the height of the tank walls. We assume that the equilibrium points satisfy the "no-spill condition" (3.2), i.e., $h^*<H_{\max}$.

Two cases have been studied in the literature (see [33, 35, 36]): the case where friction is present and surface tension is absent and the case where surface tension is present and friction is absent.

### *3.1. Preliminaries*

Define the sets $X \subset S \subset \mathbb{R}^2 \times H^1(0,L) \times L^2(0,L)$:

$$(\xi, w, h, v) \in S \Leftrightarrow \begin{cases} h \in C^0([0,L];(0,+\infty)) \cap H^1(0,L) \\ v \in C^0([0,L]) \\ \int_0^L h(x)dx = m \\ (\xi, w) \in \mathbb{R}^2, v(0) = v(L) = 0 \end{cases} \tag{3.3}$$

$$X := \left\{ (\xi, w, h, v) \in S : \max_{x \in [0,L]} \left( h(x) \right) < H_{\max} \right\} \tag{3.4}$$

The state space of system (2.9)-(2.15) is considered to be the set $X$ defined by (3.4). More specifically, we consider as state space the metric space $X \subset \mathbb{R}^2 \times H^1(0,L) \times L^2(0,L)$ with metric induced by the norm of the underlying normed linear space $\mathbb{R}^2 \times H^1(0,L) \times L^2(0,L)$, i.e.,

$$\|(\xi, w, h, v)\|_X = \left( \xi^2 + w^2 + \|h\|_2^2 + \|h'\|_2^2 + \|v\|_2^2 \right)^{1/2} \tag{3.5}$$

We define the following functionals for all $(\xi, w, h, v) \in S$:

$$V(\xi, w, h, v) := \delta E(h, v) + W(h, v) + \frac{qk^2}{2} \xi^2 + \frac{q}{2} \left( w + k\xi \right)^2 \tag{3.6}$$

$$E(h, v) := \frac{1}{2} \int_0^L h(x) v^2(x) dx + \frac{g}{2} \left\| h - h^* \chi_{[0,L]} \right\|_2^2 + \sigma \int_0^L \left( \sqrt{1 + \left( h'(x) \right)^2} - 1 \right) dx \tag{3.7}$$

$$\begin{aligned} W(h, v) &:= \frac{1}{2} \int_0^L h^{-1}(x) \left( h(x) v(x) + \mu h'(x) \right)^2 dx \\ &+ \frac{g}{2} \left\| h - h^* \chi_{[0,L]} \right\|_2^2 + \sigma \int_0^L \left( \sqrt{1 + \left( h'(x) \right)^2} - 1 \right) dx \end{aligned} \tag{3.8}$$

where $h^* = m / L$ and $\delta, k, q > 0$ are controller gains (to be selected). Moreover, we define for all $(\xi, w, h, v) \in S$ with $v \in H^1(0, L)$ the functional:

$$\tilde{V}(\xi, w, h, v) := V(\xi, w, h, v) + \left( \frac{1}{2} \|v'\|_2^2 + \gamma V(\xi, w, h, v) \right) \exp\left( \beta V(\xi, w, h, v) \right) \tag{3.9}$$

where $\beta, \gamma > 0$ are constants. We notice that:

- the functional $E$ is the mechanical energy of the liquid within the tank. Indeed, notice that $E$ is the sum of the potential energy $\frac{g}{2}\left\|h-h^{*}\chi_{[0,L]}\right\|_{2}^{2}+\sigma\int_{0}^{L}\left(\sqrt{1+\left(h'(x)\right)^{2}}-1\right)dx$ and the kinetic energy $\frac{1}{2}\int_{0}^{L}h(x)v^{2}(x)dx$ of the liquid. It should be noticed that there is no contribution to the mechanical energy of the tank-liquid interface which is a result of the boundary condition (2.14) (a constant contact angle).

- the functional $W$ is a kind of mechanical energy of the liquid within the tank and has been used extensively in the literature of isentropic or polytropic, compressible liquid flow (see [32, 37, 43, 51, 52]) as well as in [33, 34, 35, 36].

- the functional $V$ is a linear combination of the mechanical energies $E$, $W$ and a control Lyapunov function for the tank $\frac{qk^{2}}{2}\xi^{2}+\frac{q}{2}\left(w+k\xi\right)^{2}$,

- the functional $\tilde{V}$ is a functional that can provide bounds for the sup-norm of the fluid velocity (due to Agmon's inequality and the fact that $\|v'\|_{2}^{2}\le 2\tilde{V}(\xi,w,h,v)$ for all $(\xi,w,h,v)\in S$) while the other functionals cannot provide bounds for the sup-norm of the fluid velocity.

The following proposition guarantees that the sublevel sets

$$X_{V}\left(r\right):=\left\{\,(\xi,w,h,v)\in S:V(\xi,w,h,v)\le r\,\right\},\ \text{for } r\ge 0 \tag{3.10}$$

contain a neighborhood of the desired equilibrium point $\left(0,0,h^{*}\chi_{[0,L]},0\right)$ (in the topology of $X$ with metric induced by the norm $\|\ \|_{X}$ defined by (3.5)) when $r>0$.

**Proposition 1 (see Proposition 3.5 in [35] and Proposition 1 in [36]):** *Let $q,k,\delta>0$ be given and define $V$ by means of (3.6). Then for every $(\xi,w,h,v)\in S$ satisfying the inequality $\left\|(0,w,h-h^{*}\chi_{[0,L]},v)\right\|_{X}\le\varepsilon$ for some $\varepsilon>0$ with $\varepsilon<\min\left(h^{*},H_{\max}-h^{*}\right)/\sqrt{L}$, the following inequality holds:*

$$\begin{aligned}&V(\xi,w,h,v)\le\sigma(\delta+1)\sqrt{L}\left\|(\xi,w,h-h^{*}\chi_{[0,L]},v)\right\|_{X}\\&+\max\left(\mu^{2}\left(h^{*}-\varepsilon\sqrt{L}\right)^{-1},\frac{\delta+1}{2}g,\frac{\left(\delta+2\right)H_{\max}}{2},q,\frac{3qk^{2}}{2}\right)\left\|(\xi,w,h-h^{*}\chi_{[0,L]},v)\right\|_{X}^{2}\end{aligned} \tag{3.11}$$

*where $\|\cdot\|_{X}$ is defined by (3.5).*

We use the functional $V(\xi,w,h,v)$ defined by (3.6) as a CLF for system (2.9)-(2.15). However, the functional $V(\xi,w,h,v)$ can also be used for the derivation of useful bounds for the level function $h$. This is guaranteed by the following lemma.

**Lemma 1 (see Lemma 3.1 in [35] and Lemma 1 in [36]):** *Let $q,k,\delta>0$ be given. Define the increasing function $G\in C^{0}\left(\mathbb{R}\right)\cap C^{1}\left(\left(-\infty,0\right)\cup\left(0,+\infty\right)\right)$ by means of the formula*

$$G(h):=\begin{cases}\operatorname{sgn}\left(h-h^*\right)\left(\frac{2}{3}h\sqrt{h}-2h^*\sqrt{h}+\frac{4}{3}h^*\sqrt{h^*}\right) & \text{for } h>0\\ -\frac{4}{3}h^*\sqrt{h^*}+h & \text{for } h\le 0\end{cases} \tag{3.12}$$

*Let $G^{-1}:\Re\to\Re$ be the inverse function of $G$ and define*

$$Q:=1/\left(\mu\sqrt{\delta g}\right) \tag{3.13}$$

*Define $V$ by means of (18). Then for every $(\xi,w,h,v)\in S$, the following inequality holds*

$$p_1\left(V(\xi,w,h,v)\right)\le h(x)\le p_2\left(V(\xi,w,h,v)\right),\ \text{for all } x\in[0,L] \tag{3.14}$$

*where the functions $p_i:\mathbb{R}_+\to\mathbb{R}$ ($i=1,2$) are defined by the following formulae for all $s\ge 0$ when $\sigma>0$*

$$\begin{aligned}p_1(s)&:=\max\left(G^{-1}\left(-Qs\right),h^*-\sqrt{\frac{2m\left(1+\delta\right)s}{\delta\mu^2}},h^*-\sqrt{\left(\frac{s}{\sigma(\delta+1)}+L\right)^2-L^2}\right)\\ p_2(s)&:=\min\left(G^{-1}\left(Qs\right),h^*+\sqrt{\frac{2m\left(1+\delta\right)s}{\delta\mu^2}},h^*+\sqrt{\left(\frac{s}{\sigma(\delta+1)}+L\right)^2-L^2}\right)\end{aligned} \tag{3.15}$$

*and by the following formulae for all $s\ge 0$ when $\sigma=0$*

$$p_1(s):=\max\left(G^{-1}\left(-Qs\right),h^*-\sqrt{\frac{2m\left(1+\delta\right)s}{\delta\mu^2}}\right),\ p_2(s):=\min\left(G^{-1}\left(Qs\right),h^*+\sqrt{\frac{2m\left(1+\delta\right)s}{\delta\mu^2}}\right) \tag{3.16}$$

We next define an important quantity that allows us to guarantee the validity of the state constraints of the system by means of appropriate bounds on the Lyapunov functional $V$. We define

$$R:=\frac{2\mu\sqrt{\delta g h^*}}{3}\left(H_{\max}-h^*\right)\min\left(\zeta_1,\zeta_2\right) \tag{3.17}$$

where the constants $\zeta_1,\zeta_2$ are given by the following formulae

$$\zeta_1:=$$

$$\max\left(\sqrt{\frac{H_{\max}}{h^*}}-\frac{2\sqrt{h^*}}{\sqrt{H_{\max}}+\sqrt{h^*}},\frac{3\mu\sqrt{\delta}\left(H_{\max}-h^*\right)}{4m\left(1+\delta\right)\sqrt{gh^*}},\frac{3\sigma(\delta+1)\left(\sqrt{L^2+\left(H_{\max}-h^*\right)^2}-L\right)}{2\mu\sqrt{\delta gh^*}\left(H_{\max}-h^*\right)}\right)$$

$$\zeta_2:=\frac{h^*}{H_{\max}-h^*}\max\left(2,\frac{3\mu\sqrt{\delta}}{4L\sqrt{gh^*}\left(1+\delta\right)},\frac{3\sigma(\delta+1)\sqrt{h^*}}{2\mu\sqrt{\delta g}\left(\sqrt{\left(h^*\right)^2+L^2}+L\right)}\right).$$

Notice that definition (3.17), the fact that $h^*<H_{\max}$ and Lemma 1 imply that for all $(\xi,w,h,v)\in S$ with $V(\xi,w,h,v)<R$ it holds that

$$0 < p_1\left(V(\xi,w,h,v)\right) \le h(x) \le p_2\left(V(\xi,w,h,v)\right) < H_{\max}\text{, for all } x\in[0,L] \tag{3.18}$$

Therefore, the condition for not spilling out (3.2) and the positivity condition (2.15) are automatically satisfied when $(\xi,w,h,v)\in S$ with $V(\xi,w,h,v)<R$.

The following sets are also useful in what follows:

$$X_{\tilde{V}}\left(r\right):=\left\{\,(\xi,w,h,v)\in S : v\in H^1(0,L),\tilde{V}(\xi,w,h,v)\le r+\gamma r\exp(\beta r)\,\right\}\text{, for } r\ge 0 \tag{3.19}$$

Clearly, definition (3.5) implies that $X_{\tilde{V}}\left(r\right)\subseteq X_V\left(r\right)$ for all $r\ge 0$.

### 3.2. Stabilizing Controllers

The first stabilization result deals with a special case of the friction coefficient $\kappa$; namely, the case where the function $\kappa\in C^0\left((0,+\infty)\times\mathbb{R};\mathbb{R}_+\right)$ satisfies the following assumption.

**(H)** *There exists a continuous, non-increasing function* $\bar{K}:\left(0,H_{\max}\right]\to\mathbb{R}_+$ *such that for every* $\omega\in\left(0,h^*\right]$ *the following inequality holds:*

$$h^{-2}\kappa\left(h,v\right)\le\bar{K}(\omega)\text{, for all } \left(h,v\right)\in\left[\omega,H_{\max}\right]\times\mathbb{R} \tag{3.20}$$

Notice that velocity-independent friction coefficients satisfy automatically Assumption (H). For example, the velocity-independent friction coefficient $\kappa(h)=3\mu b_3/\left(3\mu+4b_3h\right)$, where $b_3>0$ is a constant, which was derived in [22], satisfies Assumption (H) with $\bar{K}(\omega)=3\omega^{-2}\mu b_3/(3\mu+4b_3\omega)$. Moreover, if there exists a constant $\kappa_{\max}>0$ such that $\kappa\left(h,v\right)\le\kappa_{\max}$ for all $\left(h,v\right)\in(0,+\infty)\times\mathbb{R}$, then the friction coefficient $\kappa(h,v)$ satisfies Assumption (H) with $\bar{K}(\omega)=\kappa_{\max}/\omega^2$. Assumption (H) allows us to achieve exponential stabilization in the set $X_V\left(r\right)$ defined by (3.16) for all $r\in\left[0,R\right)$ with $p_1(r)\ge\omega$ for any given $\omega\in\left(0,h^*\right]$.

**Theorem 1 (Special Case of Friction Coefficient and $\sigma=0$-Theorem 3.6 in [35]):** *Suppose that* $\sigma=0$ *and that the function* $\kappa\in C^0\left((0,+\infty)\times\mathbb{R};\mathbb{R}_+\right)$ *satisfies Assumption (H). Let arbitrary* $\omega\in\left(0,h^*\right]$ *be given. Pick arbitrary* $\delta>0$ *for which the following inequality holds:*

$$2g\left(\delta+1\right)>\mu\bar{K}\left(\omega\right) \tag{3.21}$$

*where* $\bar{K}:\left(0,H_{\max}\right]\to\mathbb{R}_+$ *is the function involved in Assumption (H). Pick arbitrary* $r\in[0,R)$ *for which* $p_1(r)\ge\omega$, *where* $p_1$, $R$ *are defined by (3.16), (3.17). Pick arbitrary constants* $\zeta,k,q>0$ *with*

$$k<q\,\Theta(r) \tag{3.22}$$

*where*

$$\Theta(r):=\frac{\zeta g\mu\delta\pi^2 p_1(r)}{g\mu\delta\pi^2 p_1(r)+2\zeta L\left(mgLH_{\max}(\delta+1)^2+2\mu^2\delta\pi^2 p_1(r)\right)} \tag{3.23}$$

*Define $V$ by means of (3.6) and $X_V(r)$ by means of (3.10). Then there exist constants $M,\lambda>0$ with the following property:*

**(P)** *Every solution $(\xi(t),w(t),h[t],v[t])$ of the PDE-ODE system (2.9)-(2.13) with the control law*

$$f(t)=-\zeta\left((\delta+1)\int_0^L h(t,x)v(t,x)dx+\mu\big(h(t,L)-h(t,0)\big)-q\big(w(t)+k\xi(t)\big)\right),\text{ for } t\geq 0 \qquad (3.24)$$

*that starts in $X_V(r)$, i.e., every quadruplet of functions $\xi\in C^2(\mathbb{R}_+)$, $w\in C^1(\mathbb{R}_+)$, $v\in C^0(\mathbb{R}_+\times[0,L])\cap C^1((0,+\infty)\times[0,L])$, $h\in C^1(\mathbb{R}_+\times[0,L];(0,+\infty))\cap C^2((0,+\infty)\times(0,L))$, with $(\xi(0),w(0),h[0],v[0])\in X_V(r)$ and $v[t]\in C^2((0,L))$ for each $t>0$, that satisfies (2.9)-(2.13) and (3.24) also satisfies $(\xi(t),w(t),h[t],v[t])\in X_V(r)$ for all $t\geq 0$ as well as the following estimate for all $t\geq 0$:*

$$\begin{aligned}&\left\|\left(\xi(t),w(t),h[t]-h^*\chi_{[0,L]},v[t]\right)\right\|_X\\&\leq M\exp(-\lambda t)\left\|\left(\xi(0),w(0),h[0]-h^*\chi_{[0,1]},v[0]\right)\right\|_X\end{aligned} \qquad (3.25)$$

*If, furthermore, the solution satisfies $v\in C^0\left(\mathbb{R}_+;H^1(0,L)\right)\cap C^1\left((0,+\infty);H^1(0,L)\right)$, then the following additional estimate holds for $t\geq 0$:*

$$\|v_x[t]\|_2\leq M\exp(-\lambda t)\left(\left\|(\xi(0),w(0),h[0]-h^*\chi_{[0,L]},v[0])\right\|_X+\|v_x[0]\|_2\right) \qquad (3.26)$$

Theorem 1 can also be applied to the frictionless case and we obtain the following corollary.

**Corollary 1 (The frictionless case without surface tension-Theorem 1 in [33] and Corollary 3.8 in [35]):** *Suppose that $\sigma=0$ and $\kappa\equiv 0$. Pick arbitrary $\delta>0$ and arbitrary $r\in[0,R)$, where $R>0$ is defined by (3.17). Pick arbitrary $\zeta,k,q>0$ for which (3.22) holds, where $\Theta(r)$ is defined by (3.23). Define $V$ by means of (3.6) and $X_V(r)$ by means of (3.10). Then there exist constants $M,\lambda>0$ for which property (P) holds.*

The second stabilization result deals with the general case of the friction coefficient $\kappa\in C^0\left((0,+\infty)\times\mathbb{R};\mathbb{R}_+\right)$ without surface tension. We achieve stabilization in the set $X_{\tilde{V}}(r)$ defined by (3.17) for all $r\in[0,R)$ with $p_1(r)>\omega_1$, $\sqrt{\frac{2L}{3}\left(r+\gamma r\exp(\beta r)\right)}<\omega_2$ for any given $\omega_1\in\left(0,h^*\right)$, $\omega_2>0$.

**Theorem 2 (General Case of Friction Coefficient and $\sigma=0$-Theorem 3.9 in [35]):** *Let constants $\omega_1\in\left(0,h^*\right)$, $\omega_2>0$ be given. Define*

$$\tilde{K}:=\max\left\{\,h^{-2}\kappa(h,v):\omega_1\leq h\leq H_{\max},|v|\leq\omega_2\,\right\} \qquad (3.27)$$

*Pick arbitrary constant $\delta>0$ that satisfies*

$$2g(\delta+1)>\mu\tilde{K} \qquad (3.28)$$

*Pick arbitrary constants $\zeta,q,k>0$ with*

$$k < q\tilde{\Theta} \tag{3.29}$$

*where*

$$\tilde{\Theta} := \frac{\zeta g \mu \delta \pi^2 \omega_1}{g\mu\delta\pi^2\omega_1 + 2\zeta L\left(mgLH_{\max}(\delta+1)^2 + 2\mu^2\delta\pi^2\omega_1\right)} \tag{3.30}$$

*Let* $p_1, R$ *be defined by (3.16), (3.17) and pick arbitrary constants* $\beta, \gamma > 0$ *with*

$$\gamma > 5\delta^{-1}\mu^{-1}\tilde{\alpha}^{-1}\left(H_{\max}\tilde{K}^2 + \varepsilon_1\right),\ \beta > \max\left(\frac{4\varepsilon_2}{\left(2\tilde{\alpha} + \mu\delta\gamma\omega_1\right)\omega_1^2}, \frac{20L}{3\mu^2\delta\omega_1}\right) \tag{3.31}$$

*where*

$$\varepsilon_1 := \mu^{-2}\left(\delta+1\right)g^2 H_{\max} + 3\zeta^2 L\left((\delta+1)(\delta+2)m + \delta q\right) \tag{3.32}$$

$$\varepsilon_2 := 100\delta^{-2}\mu^{-3}\left(\delta+1\right)^2 R \tag{3.33}$$

$$\tilde{\alpha} := \frac{\min\left(\mu g, 4qk^3, 4q\left(q\tilde{\Theta} - k\right), \mu\delta\right)}{2\max\left(\pi^{-2}L^2\left(\delta+2\right)\omega_1^{-1}H_{\max}, (\delta+1)gL^2 + 2\omega_1^{-1}\mu^2, qk^2, q\right)} \tag{3.34}$$

*and define* $\tilde{V}$ *by means of (3.9) and* $X_{\tilde{V}}(r)$ *by means of (3.19). Then for every* $r \in [0, R)$ *with*

$$p_1(r) > \omega_1,\ \sqrt{\frac{2L}{3}\left(r + \gamma r \exp\left(\beta r\right)\right)} < \omega_2 \tag{3.35}$$

*there exist constants* $M, \lambda > 0$ *with the following property:*

**(P')** *Every solution* $\left(\xi(t), w(t), h[t], v[t]\right)$ *of the PDE-ODE system (2.9)-(2.13) with the control law (3.24) that starts in* $X_{\tilde{V}}(r)$*, i.e., every quadruplet of functions* $\xi \in C^2\left(\mathbb{R}_+\right)$*,* $w \in C^1\left(\mathbb{R}_+\right)$*,* $v \in C^0\left(\mathbb{R}_+ \times [0,L]\right) \cap C^1\left((0,+\infty)\times[0,L]\right)$*,* $h \in C^1\left(\mathbb{R}_+ \times [0,L]; (0,+\infty)\right) \cap C^2\left((0,+\infty)\times(0,L)\right)$*, with* $\left(\xi(0), w(0), h[0], v[0]\right) \in X_{\tilde{V}}(r)$*,* $v \in C^0\left(\mathbb{R}_+; H^1\left(0,L\right)\right) \cap C^1\left((0,+\infty); H^1\left(0,L\right)\right)$ *and* $v[t] \in C^2\left((0,L)\right)$ *for each* $t > 0$*, that satisfies (2.9)-(2.13) and (3.24) also satisfies* $\left(\xi(t), w(t), h[t], v[t]\right) \in X_{\tilde{V}}(r)$ *and estimates (3.25), (3.26) for all* $t \geq 0$*.*

The third and final stabilization result deals with the frictionless case with surface tension.

**Theorem 3 (The frictionless case with surface tension-Theorem 1 in [36]):** *Suppose that* $\sigma > 0$ *and* $\kappa \equiv 0$*. Let* $p_1, R$ *be defined by (3.15), (3.17) and pick arbitrary* $r \in [0, R)$ *and arbitrary constants* $\zeta, k, q, \delta > 0$ *for which (3.22) holds, where* $\Theta(r)$ *is defined by (3.23). Define* $V$ *by means of (3.6) and* $X_V(r)$ *by means of (3.10). Then there exist constants* $M, \lambda > 0$ *with the following property:*

**(P'')** *Every solution* $\left(\xi(t), w(t), h[t], v[t]\right)$ *of the PDE-ODE system (2.9)-(2.14) with the control law (3.24) that starts in* $X_V(r)$*, i.e., every quadruplet of functions* $\xi \in C^2\left(\mathbb{R}_+\right)$*,* $w \in C^1\left(\mathbb{R}_+\right)$*,* $h \in C^1\left(\mathbb{R}_+ \times [0,L]; (0,+\infty)\right) \cap C^3\left((0,+\infty)\times(0,L)\right)$*,* $v \in C^0\left(\mathbb{R}_+ \times [0,L]\right) \cap C^1\left((0,+\infty)\times[0,L]\right)$ *with* $\left(\xi(0), w(0), h[0], v[0]\right) \in X_V(r)$*,* $v[t] \in C^2\left((0,L)\right)$ *for each* $t > 0$*, that satisfies (2.9)-(2.14) and (3.24) also satisfies* $\left(\xi(t), w(t), h[t], v[t]\right) \in X_V(r)$ *and estimate (3.25) for all* $t \geq 0$*.*

*3.3. Discussion of the Results*

While Theorem 1, Corollary 1, Theorem 2 and Theorem 3 are very different results, they all deal with the feedback law (3.24). These stabilization results show that the feedback law (3.24) guarantees exponential convergence of the tank-liquid system to the desired equilibrium point under very different situations: when surface tension is present but friction is absent and when surface tension is absent but friction may be present. Therefore, the nonlinear feedback law (3.24) is successful under very different conditions for the fluid.

It should be noticed that the feedback law (3.24) does not require the measurement of the whole liquid level and liquid velocity profile and requires the measurement of only four quantities:

- the tank position $\xi(t)$ and the tank velocity $w(t)$,
- the total liquid momentum $\int_0^L h(t,x)v(t,x)dx$, and
- the liquid level difference at the tank walls $h(t,L)-h(t,0)$.

Moreover, the feedback law (3.24) does not require the knowledge of the constant $\sigma$ (the ratio of the surface tension and liquid density) or the acceleration of gravity $g$. However, the feedback law (3.24) requires the knowledge of the kinematic viscosity of the liquid $\mu$.

On the other hand, the controller parameters $\omega,k,q,\delta>0$ are not completely free. There is one restriction for the controller parameters in all cases: (3.22) and (3.29) imply that the ratio $k/q$ must be sufficiently small. There is also an additional restriction for the controller parameters in the case where friction is present: (3.21) and (3.28) imply that $\delta>0$ must be sufficiently large.

Theorem 1, Theorem 2 and Theorem 3 are very different results because they can be applied to different notions of the solution for the closed-loop system (2.9)-(2.14) with (3.24). All theorems are applied to classical solutions. However, Theorem 1 and Corollary 1 give the additional stability estimate (3.26) when the additional regularity property $v\in C^0\left(\mathbb{R}_+;H^1(0,L)\right)\cap C^1\left((0,+\infty);H^1(0,L)\right)$ holds and Theorem 2 applies only to solutions for which the additional regularity property $v\in C^0\left(\mathbb{R}_+;H^1(0,L)\right)\cap C^1\left((0,+\infty);H^1(0,L)\right)$ holds. On the other hand, Theorem 3 does not require the regularity property $v\in C^0\left(\mathbb{R}_+;H^1(0,L)\right)\cap C^1\left((0,+\infty);H^1(0,L)\right)$ for the velocity but requires the additional regularity property $h\in C^1\left(\mathbb{R}_+\times[0,L];(0,+\infty)\right)\cap C^3\left((0,+\infty)\times(0,L)\right)$ for the liquid level.

Proposition 1 guarantees that the set $X_V\left(r\right)$ for $r>0$ contains a neighborhood of $\left(0,0,h^*\chi_{[0,L]},0\right)$ (in the topology of $X$ with metric induced by the norm $\|\ \|_X$ defined by (3.5)). Similarly, Proposition 1 and definitions (3.9), (3.19) guarantee that the set $X_{\tilde{V}}\left(r\right)$ for $r>0$ contains a neighborhood of $\left(0,0,h^*\chi_{[0,L]},0\right)$ (in the topology of $\hat{X}=\left\{\,(\xi,w,h,v)\in X: v\in H^1(0,L)\,\right\}$ with metric induced by the norm $\|(\xi,w,h,v)\|=\left(\|(\xi,w,h,v)\|_X^2+\|v'\|_2^2\right)^{1/2}$, where $\|\ \|_X$ is defined by (3.5)). Therefore, all stabilization results (Theorem 1, Corollary 1, Theorem 2 and Theorem 3) are at least local results (in the sense that they guarantee stability estimates in a neighborhood of the desired equilibrium point) with regions of attraction that contain the sets $X_V(r)$ and $X_{\tilde{V}}(r)$. The size of the sets $X_V(r)$ and $X_{\tilde{V}}(r)$ depends on $r\in\left[0,R\right)$ and on the controller parameters. However, the dependence of the size of the sets $X_V(r)$ and $X_{\tilde{V}}(r)$ on the controller parameters is very complicated and cannot be described easily. The local nature of the results means that not only do the liquid height and velocity profiles need to start near the equilibrium values, but also the tank velocity needs to start small and even the tank position needs to start near the tank set-point position.

The proofs of all the stabilization results (Theorem 1, Corollary 1, Theorem 2 and Theorem 3) and the feedback design are based on the Lyapunov functional $V(\xi,w,h,v)$ defined by (3.6), which is used as a CLF. On the other hand, the functional $\tilde{V}(\xi,w,h,v)$ defined by (3.9) is used only for the derivation of estimate (3.26) and the derivation of estimates of the sup-norm of the velocity (used in the proof of Theorem 2). This situation is something that does not happen in finite-dimensional systems: a CLF for a finite-dimensional system can provide all estimates for the solution. *Therefore, in the infinite-dimensional case we may need one functional to play the role of a CLF and a different functional for the derivation of useful stability estimates.* This situation has been observed again in [31] in the context of a nonlinear reaction-diffusion PDE, where the CLF could only provide stability estimates for the $L^2$ norm of the state and different methodologies (Stampacchia's method and an additional functional) were used for the derivation of stability estimates in the sup-norm and the $H^1$ norm of the state.

## 4. The Linearization of the Tank-Liquid System

Linearizing model (2.9)-(2.14) with $\sigma>0$ around the equilibrium point $h(x)\equiv h^*=m/L$, $v(x)\equiv 0$ and setting $\varphi=h-h^*\chi_{[0,L]}$, we obtain the following linear PDE-ODE model

$$\dot{\xi}=w \quad , \quad \dot{w}=-f \text{, for } t\geq 0 \tag{4.1}$$

$$\varphi_{tt}=c^2\varphi_{xx}-\sigma h^*\varphi_{xxxx}+\mu\varphi_{txx}-\bar{\kappa}\varphi_t \text{, for } t>0,\ x\in(0,L) \tag{4.2}$$

$$\varphi_x(t,0)=\varphi_x(t,L)=0 \text{, for } t\geq 0 \tag{4.3}$$

$$\varphi_{xxx}(t,0)=\varphi_{xxx}(t,L)=-\sigma^{-1}f(t) \text{, for } t\geq 0 \tag{4.4}$$

$$\int_0^L \varphi(t,x)dx=\int_0^L \varphi_t(t,x)dx=0 \text{, for } t\geq 0 \tag{4.5}$$

where $c=\sqrt{gh^*}$ and $\bar{\kappa}\geq 0$ are constants. The control system (4.1)-(4.5) is a system that has not been studied so far in the literature. The control input appears in the ODEs (4.1) and in the boundary condition (4.4) (boundary control). It should be noticed that the control input is a boundary input for the linearized model (4.1)-(4.5), while it is not a boundary input for the nonlinear system. This is not a result of the linearization but a result of the conversion of the linearized model to a single second-order in time PDE.

The distributed subsystem (4.2)-(4.4) is a combination of an Euler-Bernoulli beam equation and a wave equation with (internal) Kelvin-Voigt damping and possible friction. Equation (8) appears in the study of incompressible fluids flowing underground in a fractured or fissured medium where $\varphi$ is the pressure of the fluid in the porous part of the medium; see [23] on pages 217-218. The open-loop eigenvalues of the distributed subsystem (4.2)-(4.4) are the roots of the equations:

$$s^2+\left(\mu\frac{n^2\pi^2}{L^2}+\bar{\kappa}\right)s+\frac{n^2\pi^2}{L^2}\left(c^2+\sigma h^*\frac{n^2\pi^2}{L^2}\right)=0\text{, for } n=1,2,... \tag{4.6}$$

We next provide the eigenvalues in the case $\bar{\kappa}=0$. If $\mu^2\le 4\sigma h^*$ then all eigenvalues are complex and are given by the following formula for $n=1,2,...$:

$$s_n=-\frac{\mu}{2L^2}n^2\pi^2\pm i\frac{n\pi}{L}\sqrt{\frac{4\sigma h^*-\mu^2}{4L^2}n^2\pi^2+c^2} \tag{4.7}$$

If $\mu^2>4\sigma h^*$ then the eigenvalues are real for $n\ge\frac{2cL}{\pi\sqrt{\mu^2-4\sigma h^*}}$ and are given by the following formula for all $n=1,2,...$ with $n\ge\frac{2cL}{\pi\sqrt{\mu^2-4\sigma h^*}}$:

$$s_n=-\mu\frac{n^2\pi^2}{2L^2}\pm\frac{n\pi}{L}\sqrt{\frac{\mu^2-4\sigma h^*}{4L^2}n^2\pi^2-c^2} \tag{4.8}$$

In every case ($\mu^2>4\sigma h^*$ or $\mu^2\le 4\sigma h^*$), we have $\lim_{n\to+\infty}\left(\mathrm{Re}\left(s_n\right)\right)=-\infty$, indicating the strong (internal) damping that is caused by the viscosity of the fluid.

When $\bar{\kappa}>0$ then in general the eigenvalues are moved to the left in the complex plane due to the additional damping caused by the friction term $-\bar{\kappa}\varphi_t$ in the left-hand side of (4.2).

Define

$$\bar{S}=\left\{\varphi\in H^2(0,L):\varphi'(0)=\varphi'(L)=0\right\} \tag{4.9}$$

For the distributed subsystem (4.2)-(4.5) we are in a position to show the following results, which guarantee well-posedness and exponential stability.

**Theorem 4 (Existence/Uniqueness of solutions for the open-loop system):** *For every $\varphi_0\in\bar{S}\cap H^4(0,L)$, $p_0\in\bar{S}$ and $f\in C^3\left(\mathbb{R}_+\right)$ with $\varphi_0'''(0)=\varphi_0'''(L)=-\sigma^{-1}f(0)$ there exists a unique function $\varphi\in C^0\left(\mathbb{R}_+;\bar{S}\cap H^4(0,L)\right)\cap C^1\left(\mathbb{R}_+;\bar{S}\right)\cap C^2\left(\mathbb{R}_+;L^2(0,L)\right)$ with $\varphi[0]=\varphi_0$, $\varphi_t[0]=p_0$ that satisfies (4.2), (4.3) and (4.4) for all $t\ge 0$. Moreover, if $\int_0^L\varphi_0(x)dx=\int_0^L p_0(x)dx=0$ then (4.5) holds.*

**Theorem 5 (Stability properties of the open-loop system):** *There exist constants $\bar{M},\bar{\lambda},\Gamma>0$ such that for every $f\in C^0\left(\mathbb{R}_+\right)$ and for every function $\varphi\in C^0\left(\mathbb{R}_+;\bar{S}\cap H^4(0,L)\right)\cap C^1\left(\mathbb{R}_+;\bar{S}\right)\cap C^2\left(\mathbb{R}_+;L^2(0,L)\right)$ that satisfies (4.2), (4.3), (4.4) and (4.5) for all $t\ge 0$ the following estimate holds:*

$$P(t)\le\bar{M}\exp\left(-\bar{\lambda}\,t\right)P(0)+\Gamma\max_{0\le s\le t}\left(\exp\left(-\bar{\lambda}\,(t-s)\right)\left|f(s)\right|\right)\text{, for } t\ge 0 \tag{4.10}$$

*where*

$$P(t):=\left(\left\|\varphi[t]\right\|_2^2+\left\|\varphi_x[t]\right\|_2^2+\left\|\varphi_{xx}[t]\right\|_2^2+\left\|\varphi_t[t]\right\|_2^2\right)^{1/2}\text{, for } t\ge 0 \tag{4.11}$$

**Remark:** Estimate (4.11) shows that the distributed subsystem (4.2)-(4.5) satisfies the Input-to-State Stability property (see [30]) with respect to the input $f$ in the $H^2(0,L)\times L^2(0,L)$ norm of the state $(\varphi,\varphi_t)$.

## 5. Stabilization of the Linearized Model

Theorem 2 shows that the linearized model (4.1)-(4.5) is the interconnection of a double integrator ODE subsystem (recall (4.1)) with the exponentially stable distributed subsystem (4.2)-(4.5). The only connection of the ODE subsystem with the PDE subsystem is the control input that appears in both subsystems (otherwise the subsystems are completely independent). This structural feature allows us to consider two different ways of stabilizing the equilibrium point $(\xi,w)=0\in\mathbb{R}^2$, $(\varphi,\varphi_t)=0\in\bar{S}\times L^2(0,L)$ of system (4.1)-(4.5).

1st Way of Stabilization

We can stabilize exponentially the equilibrium point $(\xi,w)=0\in\mathbb{R}^2$, $(\varphi,\varphi_t)=0\in\bar{S}\times L^2(0,L)$ of system (4.1)-(4.5) by completely ignoring the liquid dynamics and using the feedback law

$$f(t)=k_1\xi(t)+k_2w(t) \tag{5.1}$$

where $k_1,k_2>0$ are constants. Then using Theorem 1 and Theorem 2 we conclude that for every $(\xi_0,w_0)\in\mathbb{R}^2$, $\varphi_0\in\bar{S}\cap H^4(0,L)$, $p_0\in\bar{S}$ with $\int_0^L\varphi_0(x)dx=\int_0^L p_0(x)dx=0$ and $\varphi_0'''(0)=\varphi_0'''(L)=-\sigma^{-1}\left(k_1\xi_0+k_2w_0\right)$ there exist unique functions $\varphi\in C^0\left(\mathbb{R}_+;\bar{S}\cap H^4(0,L)\right)\cap C^1\left(\mathbb{R}_+;\bar{S}\right)\cap C^2\left(\mathbb{R}_+;L^2(0,L)\right)$, $(\xi,w)\in C^\infty\left(\mathbb{R}_+;\mathbb{R}^2\right)$ with $(\xi(0),w(0))=(\xi_0,w_0)$, $\varphi[0]=\varphi_0$, $\varphi_t[0]=p_0$ that satisfy (4.1)-(4.5) and (5.1) for all $t\geq 0$. Moreover, there exist constants $\tilde{M},\tilde{\lambda}>0$ such that the following estimate holds:

$$\sqrt{\xi^2(t)+w^2(t)+P^2(t)}\leq\tilde{M}\exp\left(-\tilde{\lambda}t\right)\sqrt{\xi^2(0)+w^2(0)+P^2(0)}\text{, for } t\geq 0 \tag{5.2}$$

where $P(t)$ is defined by (4.11).

However, there is a problem for the feedback law (5.1). Since we have ignored completely the liquid dynamics, it is possible that the overshoot for the state component $\varphi$ (i.e., the deviation of the liquid level from the equilibrium level) is large even if the liquid starts from an almost slosh-free initial condition. In other words, *it is possible that the feedback law (5,1) agitates strongly the liquid causing sloshing during a transient period.*

2nd Way of Stabilization

The following result plays a fundamental role in what follows.

**Theorem 6 (Well-posedness of the closed-loop system under an arbitrary feedback law):** *Let $B,C\in\mathbb{R}$ be constants and let $\tilde{r},\tilde{p}\in C^0\left([0,L]\right)$ be given functions. Suppose that*

$$\langle \tilde{g},\tilde{r}\rangle \geq \frac{12L^4}{\pi^4}\sqrt{\frac{2}{L}}\sum_{n\,odd}\frac{\left|\langle \phi_n,\tilde{r}\rangle\right|}{n^4}$$
$$\langle \tilde{g}',\tilde{p}\rangle \geq \frac{12L^4}{\pi^4}\sqrt{\frac{2}{L}}\sum_{n\,odd}\frac{\left|\langle \phi_n',\tilde{p}\rangle\right|}{n^4} \tag{5.3}$$

*where* $\phi_n(x)=\sqrt{\frac{2}{L}}\cos\left(n\pi\frac{x}{L}\right)$ *for* $n=1,2,...$ *and*

$$\tilde{g}(x):=x^3-\frac{3L}{2}x^2+\frac{L^3}{4}\text{, for } x\in[0,L] \tag{5.4}$$

*Then for every* $(\xi_0,w_0)\in\mathbb{R}^2$, $\varphi_0\in\bar{S}\cap H^4(0,L)$, $u_0\in\bar{S}$ *with* $\int_0^L\varphi_0(x)dx=\int_0^L u_0(x)dx=0$, $\varphi_0'''(0)=\varphi_0'''(L)=-\sigma^{-1}\left(B\xi_0+Cw_0+\langle u_0,\tilde{r}\rangle+\langle\varphi_0',\tilde{p}\rangle\right)$ *there exist unique functions* $\varphi\in C^0\left(\mathbb{R}_+;\bar{S}\cap H^4(0,L)\right)\cap C^1\left(\mathbb{R}_+;\bar{S}\right)\cap C^2\left(\mathbb{R}_+;L^2(0,L)\right)$, $(\xi,w)\in C^1\left(\mathbb{R}_+;\mathbb{R}^2\right)$ *with* $(\xi(0),w(0))=(\xi_0,w_0)$, $\varphi[0]=\varphi_0$, $\varphi_t[0]=u_0$ *that satisfy for all* $t\geq 0$ *the equations (4.1)-(4.5) under the control law:*

$$f(t)=B\xi(t)+Cw(t)+\langle\varphi_t[t],\tilde{r}\rangle+\langle\varphi_x[t],\tilde{p}\rangle \tag{5.5}$$

**Remark:** Theorem 6 shows that the closed-loop system (4.1)-(4.5) under a feedback law of the form (5.5) is well-posed when inequalities (5.3) are valid.

We next consider the family of feedback laws given by the formula:

$$f(t)=K\left(k_5^2w(t)+k_5\xi(t)-h^*(k_3+k_4)\int_0^L x\varphi_t(t,x)dx-k_3\mu h^*\left(\varphi(t,L)-\varphi(t,0)\right)\right) \tag{5.6}$$

where $K,k_3,k_4,k_5>0$ are the control parameters with

$$k_5^{-3}<\min\left(\frac{c^2}{4k_3\mu\left(h^*\right)^2L},\frac{\mu\pi^2\left(\mu\pi^2+2K\left(h^*\right)^2L^3k_4\right)}{8K\left(h^*\right)^4L^6(k_4+k_3)^2},\frac{K}{4}\right) \tag{5.7}$$

Notice that the family of feedback laws (5.6) corresponds to the linearization of the nonlinear feedback laws (3.24) which were used for the nonlinear system (2.9)-(2.14). It should be noticed that the family of feedback laws (5.6), (5.7) is independent of the surface tension coefficient $\sigma>0$, the acceleration of gravity $g$ and the friction coefficient $\bar{\kappa}\geq 0$. Moreover, contrary to (5.1), the feedback law (5.6) is strongly affected by the liquid momentum and the liquid level. Consequently, it is expected that the feedback law (5.6) does not cause agitation of the liquid and tries to compensate between the two control objectives of bringing the tank to a specified position and having the liquid at rest. For the family of feedback laws (5.6), (5.7) we are in a position to prove the following result.

**Theorem 7 (Exponential Stabilization by means of Liquid-Dependent Feedback):** *Let* $K, k_3, k_4, k_5 > 0$ *be given constants for which (5.7) holds. Then for every* $(\xi_0, w_0) \in \mathbb{R}^2$, $\varphi_0 \in \bar{S} \cap H^4(0,L)$, $u_0 \in \bar{S}$ *with* $\int_0^L \varphi_0(x)dx = \int_0^L u_0(x)dx = 0$ *and*

$$\varphi_0'''(0) = \varphi_0'''(L) = -\sigma^{-1} K \left( k_5^2 w_0 + k_5 \xi_0 - h^*(k_3 + k_4) \int_0^L x u_0(x)dx - k_3 \mu h^* \left(\varphi_0(L) - \varphi_0(0)\right) \right)$$

*there exist unique functions* $\varphi \in C^0\left(\mathbb{R}_+; \bar{S} \cap H^4(0,L)\right) \cap C^1\left(\mathbb{R}_+; \bar{S}\right) \cap C^2\left(\mathbb{R}_+; L^2(0,L)\right)$, $(\xi, w) \in C^1\left(\mathbb{R}_+; \mathbb{R}^2\right)$ *with* $(\xi(0), w(0)) = (\xi_0, w_0)$, $\varphi[0] = \varphi_0$, $\varphi_t[0] = u_0$ *that satisfy (4.1)-(4.5) and (5.6) for all* $t \geq 0$. *Moreover, there exist constants* $\hat{M}, \hat{\lambda} > 0$ *such that the following estimate holds:*

$$\sqrt{\xi^2(t) + w^2(t) + P^2(t)} \leq \hat{M} \exp\left(-\hat{\lambda} t\right) \sqrt{\xi^2(0) + w^2(0) + P^2(0)}, \text{ for } t \geq 0 \tag{5.8}$$

*where* $P(t)$ *is defined by (4.11).*

The proof of Theorem 7 is based on the following Lyapunov functional

$$\begin{aligned} \tilde{W}(\xi, w, \varphi, \varphi_t) = {} & \frac{1}{2}\xi^2 + \frac{k_5^2}{2}\left(w + k_5^{-1}\xi\right)^2 + \frac{\mu}{K\left(h^*\right)^2 L}\left(\frac{1}{2}\|\varphi_t\|_2^2 + \frac{c^2}{2}\|\varphi'\|_2^2 + \frac{\sigma h^*}{2}\|\varphi''\|_2^2\right) \\ & + k_4\left(\frac{1}{2}\|\theta\|_2^2 + \frac{c^2}{2}\|\varphi\|_2^2 + \frac{\sigma h^*}{2}\|\varphi'\|_2^2\right) + k_3\left(\frac{1}{2}\|\theta - \mu\varphi'\|_2^2 + \frac{c^2 + \bar{\kappa}\mu}{2}\|\varphi\|_2^2 + \frac{\sigma h^*}{2}\|\varphi'\|_2^2\right) \end{aligned} \tag{5.9}$$

where

$$\theta(x) = \int_0^x \varphi_t(s)ds, \text{ for } x \in [0, L] \tag{5.10}$$

More specifically, we show that there exists a constant $\omega > 0$ such that the solutions of the closed-loop system (4.1)-(4.5) and (5.6) satisfy the differential inequality $\frac{d}{dt}\left(\tilde{W}\left(\xi(t), w(t), \varphi[t], \varphi_t[t]\right)\right) \leq -2\omega \tilde{W}\left(\xi(t), w(t), \varphi[t], \varphi_t[t]\right)$ for all $t \geq 0$.

There are major differences between the results for the linearized system and the results for the nonlinear system.

1) In the linearized case we provide existence/uniqueness results for the closed-loop system. In the nonlinear case we do not provide existence/uniqueness results for the corresponding closed-loop system.

2) In the linearized case we can study the situation where both friction and surface tension are present. In the nonlinear case, we cannot study the situation where both friction and surface tension are present.

3) There is a big difference in the state norm for which we achieve stabilization. In the nonlinear case we achieve stabilization in the state norm $\left(\xi^2 + w^2 + \left\|h - h^* \chi_{[0,L]}\right\|_2^2 + \left\|h_x\right\|_2^2 + \left\|v\right\|_2^2\right)^{1/2}$ (recall (3.5) and (3.25)). On the other hand, in the linearized case we achieve stabilization in the state norm $\left(\xi^2 + w^2 + \|\varphi\|_2^2 + \|\varphi_x\|_2^2 + \|\varphi_{xx}\|_2^2 + \|\varphi_t\|_2^2\right)^{1/2}$. In order to compare the two norms, we notice that the linearization of (2.9)-(2.14) gives $\varphi = h - h^* \chi_{[0,L]}$ as well as the equation $\varphi_t + h^* v_x = 0$. Consequently, the state norm $\left(\xi^2 + w^2 + \|\varphi\|_2^2 + \|\varphi_x\|_2^2 + \|\varphi_{xx}\|_2^2 + \|\varphi_t\|_2^2\right)^{1/2}$ corresponds to the norm $\left(\xi^2 + w^2 + \left\|h - h^* \chi_{[0,L]}\right\|_2^2 + \left\|h_x\right\|_2^2 + \left\|h_{xx}\right\|_2^2 + \left\|v_x\right\|_2^2\right)^{1/2}$. Therefore, the state norm for which stabilization is achieved in the linearized case is <u>*stronger*</u> than the state norm for which stabilization is achieved in the nonlinear case.

4) The Lyapunov functional used for the linearized case does not correspond exactly to the functionals used for the nonlinear case. The Lyapunov functional in the linearized case is a linear combination of four functionals: (i) the function $\frac{1}{2}\xi^2 + \frac{1}{2a^2}(w + a\xi)^2$ which is the Lyapunov function for the tank, (ii) the functional $\frac{1}{2}\|\theta\|_2^2 + \frac{c^2}{2}\|\varphi\|_2^2 + \frac{\sigma h^*}{2}\|\varphi_x\|_2^2$ which corresponds to the mechanical energy of the liquid, (iii) the functional $\frac{1}{2}\|\theta - \mu\varphi_x\|_2^2 + \frac{c^2 + \bar{\kappa}\mu}{2}\|\varphi\|_2^2 + \frac{\sigma h^*}{2}\|\varphi_x\|_2^2$ which corresponds to the modified mechanical energy of the liquid, and (iv) the functional $\frac{1}{2}\|\varphi_t\|_2^2 + \frac{c^2}{2}\|\varphi_x\|_2^2 + \frac{\sigma h^*}{2}\|\varphi_{xx}\|_2^2$ which is the energy of the liquid if one considers the liquid as a beam described by the beam-like equation (4.2). The first three functionals correspond to functionals which are also used for the construction of a Lyapunov functional in the nonlinear case. However, the last functional, the beam energy $\frac{1}{2}\|\varphi_t\|_2^2 + \frac{c^2}{2}\|\varphi_x\|_2^2 + \frac{\sigma h^*}{2}\|\varphi_{xx}\|_2^2$, has no nonlinear counterpart. This also explains the difference in the state norm for which stabilization is achieved.

5) As noticed above, the family of feedback laws (5.6) corresponds to the linearization of the nonlinear feedback laws (3.24) which were used for the nonlinear system (2.9)-(2.14). More specifically, when taking into account that the linearization of (2.9)-(2.14) gives $\varphi = h - h^* \chi_{[0,L]}$ and $\varphi_t + h^* v_x = 0$, we conclude that the feedback law (5.6) coincides with the feedback law (3.24) when the parameters are given by the following equations:

$$k_5 = \frac{1}{k} \quad k_3 = \frac{1}{h^* q k^2} \quad k_4 = \frac{\delta}{h^* q k^2} \quad K = q k^2 \zeta$$

Consequently, when the above equations hold the inequality (5.7) becomes (after some manipulations and using the fact that $c = \sqrt{g h^*}$):

$$k < \frac{q}{4}\min\left(\frac{g}{\mu L}, \frac{\left(\mu\pi^2 + 2\delta\zeta h^* L^3\right)\mu\pi^2}{2\zeta\left(h^*\right)^2 L^6(\delta+1)^2}, \zeta\right) \tag{5.11}$$

Therefore, we conclude that the feedback law (3.24) under the constraint (5.11) achieves stabilization for the linearization of the nonlinear system (2.9)-(2.14), no matter what the value of the surface tension coefficient $\sigma > 0$ is.

## 6. Proofs

**Proof of Theorem 4:** Let arbitrary $\varphi_0 \in \bar{S} \cap H^4(0,L)$, $p_0 \in \bar{S}$ and $f \in C^3(\mathbb{R}_+)$ with $\varphi_0'''(0) = \varphi_0'''(L) = -\sigma^{-1} f(0)$ be given. We perform the following transformation

$$\varphi(t,x) = u(t,x) + \bar{r}(x) f(t) \text{, for } t \ge 0 \text{, } x \in [0,L] \tag{6.1}$$

where $\bar{r} : [0,L] \to \mathbb{R}$ is a smooth function that satisfies

$$\bar{r}'(0) = \bar{r}'(L) = \int_0^L \bar{r}(x)dx = 0 \text{ and } \bar{r}'''(0) = \bar{r}'''(L) = -\sigma^{-1} \tag{6.2}$$

Then using (6.1), (6.2) and (4.2), (4.3), (4.4), we get the problem

$$u_{tt} = c^2 u_{xx} - \sigma h^* u_{xxxx} + \mu u_{txx} - \bar{\kappa} u_t + \bar{g} \text{, for } t > 0 \text{, } x \in (0,L) \tag{6.3}$$

$$u_x(t,0) = u_x(t,L) = u_{xxx}(t,0) = u_{xxx}(t,L) = 0 \text{, for } t \ge 0 \tag{6.4}$$

where

$$\bar{g}(t,x) = \left(c^2 \bar{r}''(x) - \sigma h^* \bar{r}^{(4)}(x)\right) f(t) + \left(\mu \bar{r}''(x) - \bar{\kappa} \bar{r}(x)\right) \dot{f}(t) - \bar{r}(x) \ddot{f}(t) \tag{6.5}$$

Defining

$$p = u_t, U = \begin{pmatrix} u \\ p \end{pmatrix}, F = \begin{pmatrix} 0 \\ \bar{g} \end{pmatrix} \tag{6.6}$$

we get from (6.3), (6.4) the initial-value problem:

$$\dot{U} + AU = F \tag{6.7}$$

with

$$U[0] = U_0 = \begin{pmatrix} u_0 \\ p_0 \end{pmatrix}, \ u_0(x) = \varphi_0(x) - \bar{r}(x) f(0) \text{, for } x \in [0,L] \tag{6.8}$$

where $A : D(A) \to X_1$ is the linear unbounded operator

$$A = \begin{bmatrix} 0 & -1 \\ \sigma h^* \dfrac{d^4}{d x^4} - c^2 \dfrac{d^2}{d x^2} & -\mu \dfrac{d^2}{d x^2} + \bar{\kappa} \end{bmatrix} \tag{6.9}$$

with $X_1$ being the real Hilbert space $X_1 = \bar{S} \times L^2(0,L)$ with scalar product

$$\left(U,\bar{U}\right)=\left\langle u,\bar{u}\right\rangle+c^2\left\langle u',\bar{u}'\right\rangle+\sigma h^*\left\langle u'',\bar{u}''\right\rangle+\left\langle p,\bar{p}\right\rangle,$$
$$\text{for all } U=\begin{pmatrix}u\\p\end{pmatrix}\in X_1, \bar{U}=\begin{pmatrix}\bar{u}\\\bar{p}\end{pmatrix}\in X_1 \tag{6.10}$$

and $D(A)\subset X_1$ being the linear space

$$D(A)=\left\{U=\begin{pmatrix}u\\p\end{pmatrix}\in \bar{S}^2 : u''\in \bar{S}\right\} \tag{6.11}$$

Notice that definitions (6.8), (6.11), the facts that $\varphi_0\in \bar{S}\cap H^4(0,L)$, $p_0\in\bar{S}$, $\varphi_0'''(0)=\varphi_0'''(L)=-\sigma^{-1}f(0)$ and properties (6.2) guarantee that $U_0\in D(A)$.

We next show that there exists $\bar{q}\geq 0$ such that the operator $A+\bar{q}\,I$, where $I$ is the identity operator, is a maximal monotone operator. Using (6.9), (6.10) we get for all $U=\begin{pmatrix}u\\p\end{pmatrix}\in D(A)$ and $\bar{q}\geq 0$:

$$\begin{aligned}&\left((A+\bar{q}\,I)U,U\right)=-\left\langle u,p\right\rangle-c^2\left\langle u',p'\right\rangle-\sigma h^*\left\langle u'',p''\right\rangle+\sigma h^*\left\langle p,u^{(4)}\right\rangle\\&-c^2\left\langle p,u''\right\rangle-\mu\left\langle p,p''\right\rangle+(\bar{q}+\kappa)\left\langle p,p\right\rangle+\bar{q}\left\langle u,u\right\rangle+\bar{q}c^2\left\langle u',u'\right\rangle+\bar{q}\sigma h^*\left\langle u'',u''\right\rangle\end{aligned} \tag{6.12}$$

Since $U=\begin{pmatrix}u\\p\end{pmatrix}\in D(A)$, it follows from (4.9), (6.11) that $u'(0)=u'(L)=p'(0)=p'(L)=u'''(0)=u'''(L)=0$. Consequently, integration by parts implies that $-\left\langle p,u''\right\rangle=\left\langle u',p'\right\rangle$, $-\left\langle p,p''\right\rangle=\left\langle p',p'\right\rangle$, $\left\langle p,u^{(4)}\right\rangle=\left\langle u'',p''\right\rangle$ and therefore we get from (6.12) for all $U=\begin{pmatrix}u\\p\end{pmatrix}\in D(A)$ and $\bar{q}\geq 0$:

$$\left((A+\bar{q}\,I)U,U\right)=-\left\langle u,p\right\rangle+\mu\|p'\|_2^2+(\bar{q}+\kappa)\|p\|_2^2+\bar{q}\|u\|_2^2+\bar{q}c^2\|u'\|_2^2+\bar{q}\sigma h^*\|u''\|_2^2 \tag{6.13}$$

The Cauchy-Schwarz inequality implies that $-\left\langle u,p\right\rangle\geq-\|u\|_2\|p\|_2\geq-\frac{1}{2}\|p\|_2^2-\frac{1}{2}\|u\|_2^2$. Consequently, we get from (6.13) for all $U=\begin{pmatrix}u\\p\end{pmatrix}\in D(A)$ and $\bar{q}\geq 1/2$:

$$\left((A+\bar{q}\,I)U,U\right)\geq 0 \tag{6.14}$$

Let arbitrary $\begin{pmatrix}f_1\\f_2\end{pmatrix}\in X_1$ be given. By virtue of (6.9), the equation $(A+(\bar{q}+1)I)U=\begin{pmatrix}f_1\\f_2\end{pmatrix}$ gives

$$\begin{aligned}&(\bar{q}+1)u-p=f_1\\&\sigma h^*u^{(4)}-c^2u''-\mu p''+(\bar{q}+1+\kappa)p=f_2\end{aligned} \tag{6.15}$$

The system of equations (6.15) gives the equation

$$\sigma h^* u^{(4)} - \left(c^2 + \mu(\bar{q}+1)\right)u'' + (\bar{q}+1)(\bar{q}+1+\bar{\kappa})u = f_3 \tag{6.16}$$

where $f_3 = f_2 + (\bar{q}+1+\bar{\kappa}) f_1 - \mu f_1''$. Using Fourier series we find that for every $f_3 \in L^2(0,L)$ and every $\bar{q} \geq 0$, equation (6.16) has a solution $u \in \bar{S}$ with $u'' \in \bar{S}$ which is given by the following formula for $x \in [0,L]$:

$$\begin{aligned} u(x) &= \frac{1}{(\bar{q}+1)(\bar{q}+1+\bar{\kappa})L}\int_0^L f_3(z)dz \\ &+\sum_{n=1}^{\infty}\frac{2L^3\cos\left(n\pi\frac{x}{L}\right)\int_0^L f_3(z)\cos\left(n\pi\frac{z}{L}\right)dz}{\sigma h^* n^4\pi^4 + L^2\left(c^2+\mu(\bar{q}+1)\right)n^2\pi^2 + L^4(\bar{q}+1)(\bar{q}+1+\bar{\kappa})} \end{aligned} \tag{6.17}$$

Using (6.15), (6.16) we conclude that for every $\begin{pmatrix} f_1 \\ f_2 \end{pmatrix} \in X_1$ and every $\bar{q} \geq 0$ there exists $U = \begin{pmatrix} u \\ p \end{pmatrix} \in D(A)$ such that $(A + (\bar{q}+1)I)U = \begin{pmatrix} f_1 \\ f_2 \end{pmatrix}$. Therefore, using (6.14) we conclude that the operator $A + \bar{q}\, I$ is a maximal monotone operator for $\bar{q} \geq 1/2$.

The proof is finished by applying Theorem 7.10 on page 198 in [9]. The proof is complete. $\lhd$

**Proof of Theorem 5:** Let $f \in C^0(\mathbb{R}_+)$ and an arbitrary function $\varphi \in C^0\left(\mathbb{R}_+; \bar{S} \cap H^4(0,L)\right) \cap C^1\left(\mathbb{R}_+; \bar{S}\right) \cap C^2\left(\mathbb{R}_+; L^2(0,L)\right)$ that satisfies (4.2)-(4.5) for all $t \geq 0$ be given. Define

$$\theta(t,x) = \int_0^x \varphi_t(t,s)ds \text{, for } t \geq 0,\ x \in [0,L] \tag{6.18}$$

Using (4.2)-(4.5) and definition (6.18) we conclude that the following equations hold:

$$\theta_t = c^2\varphi_x - \sigma h^*\varphi_{xxx} - h^* f + \mu\varphi_{tx} - \bar{\kappa}\,\theta \text{. for } t \geq 0 \tag{6.19}$$

$$\theta(t,0) = \theta(t,L) = 0 \text{, for } t \geq 0 \tag{6.20}$$

Using (4.2)-(4.4), (6.19), (6.20) and integration by parts we conclude that the following equations hold for $t \geq 0$:

$$\frac{d}{dt}\left(\frac{1}{2}\|\varphi_t[t]\|_2^2 + \frac{c^2}{2}\|\varphi_x[t]\|_2^2 + \frac{\sigma h^*}{2}\|\varphi_{xx}[t]\|_2^2\right) = -\mu\|\varphi_{tx}[t]\|_2^2 - \bar{\kappa}\|\varphi_t[t]\|_2^2 + h^*\left\langle \varphi_{tx}[t], \chi_{[0,L]}\right\rangle f(t) \tag{6.21}$$

$$\frac{d}{dt}\left(\frac{1}{2}\|\theta[t]\|_2^2 + \frac{c^2}{2}\|\varphi[t]\|_2^2 + \frac{\sigma h^*}{2}\|\varphi_x[t]\|_2^2\right) = -\bar{\kappa}\|\theta[t]\|_2^2 - \mu\|\varphi_t[t]\|_2^2 - h^*\left\langle \theta[t], \chi_{[0,L]}\right\rangle f(t) \tag{6.22}$$

$$\frac{d}{dt}\left(\frac{1}{2}\left\|\theta[t]-\mu\varphi_x[t]\right\|_2^2+\frac{c^2+\bar{\kappa}\mu}{2}\left\|\varphi[t]\right\|_2^2+\frac{\sigma h^*}{2}\left\|\varphi_x[t]\right\|_2^2\right)$$
$$=-\mu c^2\left\|\varphi_x[t]\right\|_2^2-\mu\sigma h^*\left\|\varphi_{xx}[t]\right\|_2^2-\bar{\kappa}\left\|\theta[t]\right\|_2^2 \tag{6.23}$$
$$-h^*\left\langle \theta[t]-\mu\varphi_x[t],\chi_{[0,L]}\right\rangle f(t)$$

Define the mapping:

$$V_1(t)=\frac{1}{2}\left\|\varphi_t[t]\right\|_2^2+\frac{c^2+2\sigma h^*}{2}\left\|\varphi_x[t]\right\|_2^2+\frac{\sigma h^*}{2}\left\|\varphi_{xx}[t]\right\|_2^2+\frac{R}{2}\left\|\theta[t]\right\|_2^2$$
$$+\frac{2c^2+\bar{\kappa}\mu}{2}\left\|\varphi[t]\right\|_2^2+\frac{1}{2}\left\|\theta[t]-\mu\varphi_x[t]\right\|_2^2 \tag{6.24}$$

We notice that $V_1(t)$ as defined by (6.24) is nothing else but the sum of the quantities whose time derivatives appear in the left-hand sides of (6.21), (6.22) and (6.23). Therefore, we get from (6.21), (6.22), (6.23) and (6.24) for all $t \ge 0$:

$$\dot{V}_1(t)=-\mu\left\|\varphi_{tx}[t]\right\|_2^2-\left(\mu+\bar{\kappa}\right)\left\|\varphi_t[t]\right\|_2^2-2\bar{\kappa}\left\|\theta[t]\right\|_2^2$$
$$-\mu c^2\left\|\varphi_x[t]\right\|_2^2-\mu\sigma h^*\left\|\varphi_{xx}[t]\right\|_2^2-h^*\left\langle 2\theta[t]-\mu\varphi_x[t]-\varphi_{tx}[t],\chi_{[0,L]}\right\rangle f(t) \tag{6.25}$$

Using Wirtinger's inequality, (4.3), (4.5) and (6.18), (6.20) we get for all $t \ge 0$:

$$\left\|\varphi_x[t]\right\|_2^2 \le \frac{L^2}{\pi^2}\left\|\varphi_{xx}[t]\right\|_2^2 \tag{6.26}$$

$$\left\|\theta[t]\right\|_2^2 \le \frac{L^2}{\pi^2}\left\|\varphi_t[t]\right\|_2^2 \tag{6.27}$$

$$\left\|\varphi[t]\right\|_2^2 \le \frac{L^2}{\pi^2}\left\|\varphi_x[t]\right\|_2^2 \tag{6.28}$$

$$\left\|\varphi_t[t]\right\|_2^2 \le \frac{L^2}{\pi^2}\left\|\varphi_{tx}[t]\right\|_2^2 \tag{6.29}$$

Using definitions (4.11), (6.24), (6.18) and (6.27) we conclude that there exist constants $K_2 > K_1 > 0$ (independent of $t \ge 0$ and the solution $\varphi$) such that the following inequalities hold for all $t \ge 0$:

$$K_1 P^2(t) \le V_1(t) \le K_2 P^2(t) \tag{6.30}$$

The Cauchy-Schwarz inequality and (6.27) gives the inequalities $\left|\left\langle \varphi_{tx}[t],\chi_{[0,L]}\right\rangle\right| \le \sqrt{L}\left\|\varphi_{tx}[t]\right\|_2$, $\left|\left\langle \varphi_x[t],\chi_{[0,L]}\right\rangle\right| \le \sqrt{L}\left\|\varphi_x[t]\right\|_2$ and $\left|\left\langle \theta[t],\chi_{[0,L]}\right\rangle\right| \le \frac{L\sqrt{L}}{\pi}\left\|\varphi_t[t]\right\|_2$ for all $t \ge 0$. Consequently, using the previous inequalities and the following elementary inequalities

$$2h^*\frac{L\sqrt{L}}{\pi}\left\|\varphi_t[t]\right\|_2 \left|f(t)\right| \le \frac{\mu+\bar{\kappa}}{2}\left\|\varphi_t[t]\right\|_2^2+\frac{2\left(h^*\right)^2 L^3}{\pi^2\left(\mu+\bar{\kappa}\right)}\left|f(t)\right|^2$$

$$h^*\mu\sqrt{L}\,\|\varphi_x[t]\|_2\,|f(t)| \le \frac{\mu c^2}{2}\|\varphi_x[t]\|_2^2 + \frac{\left(h^*\right)^2\mu L}{2c^2}|f(t)|^2$$

$$h^*\sqrt{L}\,\|\varphi_{tx}[t]\|_2\,|f(t)| \le \frac{\mu}{2}\|\varphi_{tx}[t]\|_2^2 + \frac{\left(h^*\right)^2 L}{2\mu}|f(t)|^2$$

we obtain from (6.25) for all $t \ge 0$:

$$\begin{aligned}
\dot{V}_1(t) &\le -\mu\|\varphi_{tx}[t]\|_2^2 - \left(\mu+\bar{\kappa}\right)\|\varphi_t[t]\|_2^2 - 2\bar{\kappa}\|\theta[t]\|_2^2 \\
&-\mu c^2\|\varphi_x[t]\|_2^2 - \mu\sigma h^*\|\varphi_{xx}[t]\|_2^2 + 2h^*\frac{L\sqrt{L}}{\pi}\|\varphi_t[t]\|_2\,|f(t)| \\
&+h^*\mu\sqrt{L}\,\|\varphi_x[t]\|_2\,|f(t)| + h^*\sqrt{L}\,\|\varphi_{tx}[t]\|_2\,|f(t)| \\
&\le -\frac{\mu}{2}\|\varphi_{tx}[t]\|_2^2 - \frac{\mu+\bar{\kappa}}{2}\|\varphi_t[t]\|_2^2 - \frac{\mu c^2}{2}\|\varphi_x[t]\|_2^2 \\
&-\mu\sigma h^*\|\varphi_{xx}[t]\|_2^2 + \frac{\left(h^*\right)^2 L}{2\mu}\left(\frac{4\mu L^2}{\pi^2\left(\mu+\bar{\kappa}\right)} + \frac{\mu^2}{c^2} + 1\right)|f(t)|^2
\end{aligned} \tag{6.31}$$

Using (6.28), (6.29) and (6.31) we get for all $t \ge 0$:

$$\begin{aligned}
\dot{V}_1(t) &\le -\frac{\mu\pi^2 + (\mu+\bar{\kappa})L^2}{2L^2}\|\varphi_t[t]\|_2^2 - \frac{\mu c^2}{4}\|\varphi_x[t]\|_2^2 - \frac{\mu c^2}{4}\|\varphi_x[t]\|_2^2 \\
&-\mu\sigma h^*\|\varphi_{xx}[t]\|_2^2 + \frac{\left(h^*\right)^2 L}{2\mu}\left(\frac{4\mu L^2}{\pi^2\left(\mu+\bar{\kappa}\right)} + \frac{\mu^2}{c^2} + 1\right)|f(t)|^2 \\
&\le -\frac{\mu\pi^2 + (\mu+\bar{\kappa})L^2}{2L^2}\|\varphi_t[t]\|_2^2 - \frac{\mu c^2\pi^2}{4L^2}\|\varphi[t]\|_2^2 - \frac{\mu c^2}{4}\|\varphi_x[t]\|_2^2 \\
&-\mu\sigma h^*\|\varphi_{xx}[t]\|_2^2 + \frac{\left(h^*\right)^2 L}{2\mu}\left(\frac{4\mu L^2}{\pi^2\left(\mu+\bar{\kappa}\right)} + \frac{\mu^2}{c^2} + 1\right)|f(t)|^2
\end{aligned} \tag{6.32}$$

Definition (4.11) and (6.32) implies that there exists a constant $K_3 > 0$ (independent of $t \ge 0$ and the solution $\varphi$) such that the following inequality holds for all $t \ge 0$:

$$\dot{V}(t) \le -K_3 P^2(t) + \frac{\left(h^*\right)^2 L}{2\mu}\left(\frac{4\mu L^2}{\pi^2\left(\mu+\bar{\kappa}\right)} + \frac{\mu^2}{c^2} + 1\right)|f(t)|^2 \tag{6.33}$$

Using (6.30) and (6.33) we conclude that there exists a constant $K_4 > 0$ (independent of $t \ge 0$ and the solution $\varphi$) such that the following inequality holds for all $t \ge 0$:

$$\dot{V}(t) \le -K_4 V(t) + \frac{\left(h^*\right)^2 L}{2\mu}\left(\frac{4\mu L^2}{\pi^2\left(\mu+\bar{\kappa}\right)} + \frac{\mu^2}{c^2} + 1\right)|f(t)|^2 \tag{6.34}$$

Estimate (4.10) for appropriate constants $\bar{M}, \bar{\lambda}, \Gamma > 0$ is a direct consequence of differential inequality (6.34) and inequalities (6.30). The proof is complete. ◁

**Proof of Theorem 6:** The proof follows a similar notation with the proof of Theorem 4 (for example, we have states $u, U$, scalar product $(\bullet,\bullet)$, identity operator $I$, etc.). However, the reader should not be tempted, by an overlapping notation for different quantities, to compare the different quantities in the proofs. The proofs of Theorem 4 and Theorem 6 are completely independent.

Define the real Hilbert space

$$X_2 = \left\{ (\xi, w, \varphi, u) : (\xi, w) \in \mathbb{R}^2, \varphi \in \bar{S}, u \in L^2(0,L), \int_0^L \varphi(x)dx = \int_0^L u(x)dx = 0 \right\} \tag{6.35}$$

with scalar product

$$\left(U, \bar{U}\right) = \xi\bar{\xi} + w\bar{w} + \langle \varphi, \bar{\varphi} \rangle + c^2 \langle \varphi', \bar{\varphi}' \rangle + \sigma h^* \langle \varphi'', \bar{\varphi}'' \rangle + \langle u, \bar{u} \rangle,$$
$$\text{for all } U = (\xi, w, \varphi, u) \in X_2, \bar{U} = (\bar{\xi}, \bar{w}, \bar{\varphi}, \bar{u}) \in X_2 \tag{6.36}$$

Define the linear unbounded operator $\tilde{A} : D(\tilde{A}) \to X_2$ by means of the formula

$$\tilde{A}U = \begin{pmatrix} -w \\ B\xi + Cw + \langle u, \tilde{r} \rangle + \langle \varphi', \tilde{p} \rangle \\ -u \\ -c^2\varphi'' + \sigma h^* \varphi^{(4)} - \mu u'' + \kappa u \end{pmatrix}, \text{ for all } U = (\xi, w, \varphi, u) \in D(\tilde{A}) \tag{6.37}$$

where $D(\tilde{A}) \subset X_2$ is the linear subspace

$$D(\tilde{A}) = \left\{ (\xi, w, \varphi, u) \in X_2 : \begin{array}{c} \varphi \in H^4(0,L), u \in \bar{S} \\ \sigma\varphi'''(0) = \sigma\varphi'''(L) = -B\xi - Cw - \langle u, \tilde{r} \rangle - \langle \varphi', \tilde{p} \rangle \end{array} \right\} \tag{6.38}$$

It is clear from definitions (6.37), (6.38) that we are seeking for a solution of the initial-value problem

$$\dot{U} + \tilde{A}U = 0 \tag{6.39}$$

with

$$U[0] = (\xi_0, w_0, \varphi_0, u_0) \tag{6.40}$$

Notice that definitions (6.35), (6.38), (6.40) and the facts that $(\xi_0, w_0) \in \mathbb{R}^2$, $\varphi_0 \in \bar{S} \cap H^4(0,L)$, $u_0 \in \bar{S}$ with $\int_0^L \varphi_0(x)dx = \int_0^L u_0(x)dx = 0$, $\varphi_0'''(0) = \varphi_0'''(L) = -\sigma^{-1}\left(B\xi_0 + Cw_0 + \langle u_0, \tilde{r} \rangle + \langle \varphi_0', \tilde{p} \rangle\right)$ imply that $U[0] = (\xi_0, w_0, \varphi_0, u_0) \in D(\tilde{A})$.

The theorem is proved by applying the Hille-Yosida Theorem (Theorem 7.4 on page 185 in [9] and Remark 6 on page 190 in [9]) to the initial-value problem (6.39), (6.40). To this purpose, it suffices to show that the linear unbounded operator $\tilde{A}+\bar{q}\,I$ where $\tilde{A}:D(\tilde{A})\to X_2$ is defined by (6.37), (6.38), is a maximal monotone operator for some $\bar{q}\ge 0$.

Indeed, by virtue of (4.9), (6.35), (6.36), (6.37), (6.38) and by using integration by parts, the Cauchy-Schwarz inequality and the fact that $\left|u(L)-u(0)\right|\le\sqrt{L}\,\|u'\|_2$ (a consequence of the Cauchy-Schwarz inequality and the fact that $u(L)-u(0)=\int_0^L u'(x)dx$), we have for all $U=(\xi,w,\varphi,u)\in D(\tilde{A})$:

$$\begin{aligned}
&\left(\tilde{A}U,U\right)=-\xi w+B\xi w+Cw^2+w\langle u,\tilde{r}\rangle+w\langle\varphi',\tilde{p}\rangle\\
&-\langle\varphi,u\rangle-c^2\langle\varphi',u'\rangle-\sigma h^*\langle\varphi'',u''\rangle-c^2\langle u,\varphi''\rangle\\
&+\sigma h^*\left\langle u,\varphi^{(4)}\right\rangle-\mu\langle u,u''\rangle+\bar{\kappa}\|u\|_2^2\\
&=(B-1)\xi w+Cw^2+w\langle u,\tilde{r}\rangle+w\langle\varphi',\tilde{p}\rangle+\mu\|u'\|_2^2+\bar{\kappa}\|u\|_2^2\\
&-\langle\varphi,u\rangle-h^*\left(u(L)-u(0)\right)\left(B\xi+Cw+\langle u,\tilde{r}\rangle+\langle\varphi',\tilde{p}\rangle\right)\\
&\ge-|B-1||\xi||w|-|C|w^2-\|\tilde{r}\|_2|w|\|u\|_2-\|\tilde{p}\|_2|w|\|\varphi'\|_2+\mu\|u'\|_2^2+\bar{\kappa}\|u\|_2^2\\
&-\|\varphi\|_2\|u\|_2-|B|h^*\sqrt{L}\|u'\|_2|\xi|\\
&-h^*|C|\sqrt{L}\|u'\|_2|w|-h^*\sqrt{L}\|u'\|_2\|\tilde{r}\|_2\|u\|_2-h^*\sqrt{L}\|u'\|_2\|\tilde{p}\|_2\|\varphi'\|_2
\end{aligned}\qquad(6.41)$$

Using the inequalities

$$\begin{aligned}
&|B|h^*\sqrt{L}\|u'\|_2|\xi|\le\frac{\mu}{4}\|u'\|_2^2+\frac{1}{\mu}B^2\left(h^*\right)^2L\xi^2\\
&h^*|C|\sqrt{L}\|u'\|_2|w|\le\frac{\mu}{4}\|u'\|_2^2+\frac{1}{\mu}C^2\left(h^*\right)^2Lw^2\\
&h^*\sqrt{L}\|u'\|_2\|\tilde{r}\|_2\|u\|_2\le\frac{\mu}{4}\|u'\|_2^2+\frac{1}{\mu}\left(h^*\right)^2L\|\tilde{r}\|_2^2\|u\|_2^2\\
&h^*\sqrt{L}\|u'\|_2\|\tilde{p}\|_2\|\varphi'\|_2\le\frac{\mu}{4}\|u'\|_2^2+\frac{1}{\mu}\left(h^*\right)^2L\|\tilde{p}\|_2^2\|\varphi'\|_2^2
\end{aligned}$$

$$\begin{aligned}
&|B-1||\xi||w|\le\frac{|B-1|}{2}\xi^2+\frac{|B-1|}{2}w^2\\
&\|\tilde{r}\|_2|w|\|u\|_2\le\frac{\|\tilde{r}\|_2}{2}\|u\|_2^2+\frac{\|\tilde{r}\|_2}{2}w^2\\
&\|\tilde{p}\|_2|w|\|\varphi'\|_2\le\frac{\|\tilde{p}\|_2}{2}\|\varphi'\|_2^2+\frac{\|\tilde{p}\|_2}{2}w^2\\
&\|\varphi\|_2\|u\|_2\le\frac{1}{2}\|u\|_2^2+\frac{1}{2}\|\varphi\|_2^2
\end{aligned}$$

we get from (6.36) and (6.41):

$$\left((\tilde{A}+\bar{q}\,I)U,U\right)\geq\left(\bar{q}-\frac{|B-1|}{2}-\frac{L}{\mu}B^2\left(h^*\right)^2\right)\xi^2+\bar{q}\sigma h^*\left\|\varphi''\right\|_2^2$$
$$+\left(\bar{q}-\frac{|B-1|+\left\|\tilde{r}\right\|_2+\left\|\tilde{p}\right\|_2}{2}-|C|-\frac{L}{\mu}C^2\left(h^*\right)^2\right)w^2+\left(\bar{q}-\frac{1}{2}\right)\left\|\varphi\right\|_2^2 \tag{6.42}$$
$$+\left(\bar{q}c^2-\frac{\left\|\tilde{p}\right\|_2}{2}-\frac{L}{\mu}\left(h^*\right)^2\left\|\tilde{p}\right\|_2^2\right)\left\|\varphi'\right\|_2^2+\left(\bar{q}-\frac{\left\|\tilde{r}\right\|_2+1}{2}-\frac{L}{\mu}\left(h^*\right)^2\left\|\tilde{r}\right\|_2^2\right)\left\|u\right\|_2^2$$

It follows from (6.42) that for $\bar{q}\geq 0$ sufficiently large it holds that

$$\left((\tilde{A}+\bar{q}\,I)U,U\right)\geq 0 \tag{6.43}$$

Let arbitrary $f_1,f_2\in\mathbb{R}$, $f_3\in H^2(0,L), f_4\in L^2(0,L)$ with $f_3'(0)=f_3'(L)=0$ and $\int_0^L f_3(x)dx=\int_0^L f_4(x)dx=0$ be given. We investigate the existence of $U=\left(\xi,w,\varphi,u\right)\in D(\tilde{A})$ with $\tilde{A}U+(\bar{q}+1)U=\left(f_1,f_2,f_3,f_4\right)$. Using (6.37) we get the equations:

$$\begin{aligned}
&(\bar{q}+1)\xi-w=f_1\\
&(\bar{q}+1)w+B\xi+Cw+\left\langle u,\tilde{r}\right\rangle+\left\langle \varphi',\tilde{p}\right\rangle=f_2\\
&(\bar{q}+1)\varphi-u=f_3\\
&(\bar{q}+1)u-c^2\varphi''+\sigma h^*\varphi^{(4)}-\mu u''+\kappa u=f_4
\end{aligned} \tag{6.44}$$

For $\bar{q}\geq 0$ sufficiently large (so that $\left(\bar{q}+1+C\right)(\bar{q}+1)+B>0$ ) we get from (6.44):

$$\begin{aligned}
&\xi=\frac{\bar{q}+1+C}{s(\bar{q})}f_1+\frac{1}{s(\bar{q})}f_2+\frac{1}{s(\bar{q})}\left\langle f_3,\tilde{r}\right\rangle-\frac{\bar{q}+1}{s(\bar{q})}\left\langle\varphi,\tilde{r}\right\rangle-\frac{1}{s(\bar{q})}\left\langle\varphi',\tilde{p}\right\rangle\\
&w=-\frac{B}{s(\bar{q})}f_1+\frac{\bar{q}+1}{s(\bar{q})}f_2+\frac{\bar{q}+1}{s(\bar{q})}\left\langle f_3,\tilde{r}\right\rangle-\frac{(\bar{q}+1)^2}{s(\bar{q})}\left\langle\varphi,\tilde{r}\right\rangle-\frac{\bar{q}+1}{s(\bar{q})}\left\langle\varphi',\tilde{p}\right\rangle\\
&u=(\bar{q}+1)\varphi-f_3
\end{aligned} \tag{6.45}$$

where $s(\bar{q})=\left(\bar{q}+1+C\right)(\bar{q}+1)+B$ and $\varphi\in H^4(0,L)$ is a function that satisfies

$$\sigma h^*\varphi^{(4)}-\left(c^2+\mu(\bar{q}+1)\right)\varphi''+(\bar{q}+1)(\bar{q}+1+\kappa)\varphi=f_5 \tag{6.46}$$

$$\begin{aligned}
&\varphi'(0)=\varphi'(L)=0\\
&\varphi'''(0)=\varphi'''(L)=Z-\frac{(\bar{q}+1)^3}{\sigma s(\bar{q})}\left\langle\varphi,\tilde{r}\right\rangle-\frac{(\bar{q}+1)^2}{\sigma s(\bar{q})}\left\langle\varphi',\tilde{p}\right\rangle
\end{aligned} \tag{6.47}$$

with $Z=-\frac{(\bar{q}+1)B}{\sigma s(\bar{q})}f_1-\frac{C(\bar{q}+1)+B}{\sigma s(\bar{q})}f_2+\frac{(\bar{q}+1)^2}{\sigma s(\bar{q})}\langle f_3,\tilde{r}\rangle$ and $f_5=f_4-\mu f_3''+(\bar{q}+1+\kappa)f_3$. Notice that by virtue of (6.46), (6.47) and the facts that $f_3'(0)=f_3'(L)=0$ and $\int_0^L f_3(x)dx=\int_0^L f_4(x)dx=0$, it follows that $\int_0^L \varphi(x)dx=0$. Therefore, a solution of the boundary-value problem (6.46), (6.47) gives (by means of (6.45)) a solution $U=(\xi,w,\varphi,u)\in D(\tilde{A})$ of the equation $\tilde{A}U+(\bar{q}+1)U=(f_1,f_2,f_3,f_4)$.

Therefore, we next finish the proof by showing that for sufficiently large $\bar{q}\geq 0$ (so that $s(\bar{q})=(\bar{q}+1+C)(\bar{q}+1)+B>0$), the boundary-value problem (6.46), (6.47) has a solution $\varphi\in H^4(0,L)$ for every $Z\in\mathbb{R}$ and every $f_5\in L^2(0,L)$ with $\int_0^L f_5(x)dx=0$. We look for a solution of the form $\varphi(x)=a_0\left(x^3-\frac{3L}{2}x^2+\frac{L^3}{4}\right)+\sum_{n=1}^{\infty}a_n\phi_n(x)$, where $a_n\in\mathbb{R}$ for $n=1,2,...$ and $\phi_n(x)=\sqrt{\frac{2}{L}}\cos\left(n\pi\frac{x}{L}\right)$ for $n=1,2,....$ Substituting this expression in (6.46), (6.47) and recalling (5.4) we obtain:

$$a_n=\frac{L^4\langle f_5,\phi_n\rangle}{\tilde{\theta}(n,\bar{q})}+6a_0\sqrt{\frac{2}{L}}\frac{L^6\left((-1)^n-1\right)}{\tilde{\theta}(n,\bar{q})}b(n,\bar{q})\text{, for } n=1,2,... \tag{6.48}$$

$$\left(6+\frac{(\bar{q}+1)^3}{\sigma s(\bar{q})}\tilde{\Gamma}+\frac{(\bar{q}+1)^2}{\sigma s(\bar{q})}\tilde{\Delta}\right)a_0=Z-\frac{(\bar{q}+1)^2L^4}{\sigma s(\bar{q})}\sum_{n=1}^{\infty}\frac{\langle f_5,\phi_n\rangle}{\tilde{\theta}(n,\bar{q})}\left((\bar{q}+1)\langle\phi_n,\tilde{r}\rangle+\langle\phi_n',\tilde{p}\rangle\right) \tag{6.49}$$

where

$$\tilde{\theta}(n,\bar{q})=\sigma h^* n^4\pi^4+L^2\left(c^2+\mu(\bar{q}+1)\right)n^2\pi^2+(\bar{q}+1)(\bar{q}+1+\kappa)L^4 \tag{6.50}$$

$$b(n,\bar{q})=\frac{\left(c^2+\mu(\bar{q}+1)\right)n^2\pi^2+L^2(\bar{q}+1)(\bar{q}+1+\kappa)}{n^4\pi^4} \tag{6.51}$$

$$\begin{aligned}\tilde{\Gamma}&=\langle\tilde{g},\tilde{r}\rangle+6\sqrt{\frac{2}{L}}\sum_{n=1}^{\infty}\frac{L^6\left((-1)^n-1\right)}{\tilde{\theta}(n,\bar{q})}b(n,\bar{q})\langle\phi_n,\tilde{r}\rangle\\ \tilde{\Delta}&=\langle\tilde{g}',\tilde{p}\rangle+6\sqrt{\frac{2}{L}}\sum_{n=1}^{\infty}\frac{L^6\left((-1)^n-1\right)}{\tilde{\theta}(n,\bar{q})}b(n,\bar{q})\langle\phi_n',\tilde{p}\rangle\end{aligned} \tag{6.52}$$

Notice that the fact (6.48), (6.50), (6.51) and the fact that $f_5\in L^2(0,L)$ guarantee that $\varphi(x)=a_0\left(x^3-\frac{3L}{2}x^2+\frac{L^3}{4}\right)+\sum_{n=1}^{\infty}a_n\phi_n(x)$ is indeed a function of class $H^4(0,L)$. Therefore, the solvability of the boundary-value problem (6.46), (6.47) depends on the solvability of equation (6.49). Since $\tilde{\theta}(n,\bar{q}) = n^4\pi^4\left(\sigma h^* + L^2 b(n,\bar{q})\right)$ (a consequence of (6.50) and (6.51)) and since $b(n,\bar{q}) \geq 0$, we get from (6.52):

$$\begin{aligned} \tilde{\Gamma} &\geq \langle \tilde{g}, \tilde{r} \rangle - \frac{12L^4}{\pi^4}\sqrt{\frac{2}{L}} \sum_{n\ odd} \frac{\left|\langle \phi_n, \tilde{r} \rangle\right|}{n^4} \\ \tilde{\Delta} &\geq \langle \tilde{g}', \tilde{p} \rangle - \frac{12L^4}{\pi^4}\sqrt{\frac{2}{L}} \sum_{n\ odd} \frac{\left|\langle \phi_n', \tilde{p} \rangle\right|}{n^4} \end{aligned} \tag{6.53}$$

Inequalities (5.3) and (6.53) guarantee that $\tilde{\Gamma} \geq 0$ and $\tilde{\Delta} \geq 0$. Consequently, equation (6.49) is solvable.

Thus, we conclude that for sufficiently large $\bar{q} \geq 0$ (so that $s(\bar{q}) = \left(\bar{q}+1+C\right)(\bar{q}+1) + B > 0$) and for every $f_1, f_2 \in \mathbb{R}$, $f_3 \in H^2(0,L), f_4 \in L^2(0,L)$ with $f_3'(0) = f_3'(L) = 0$ and $\int_0^L f_3(x)dx = \int_0^L f_4(x)dx = 0$ there exists $U = \left(\xi, w, \varphi, u\right) \in D(\tilde{A})$ with $\tilde{A}U + (\bar{q}+1)U = \left(f_1, f_2, f_3, f_4\right)$.

The proof is complete. ◁

**Proof of Theorem 7:** Let $(\xi_0, w_0) \in \mathbb{R}^2$, $\varphi_0 \in \bar{S} \cap H^4(0,L)$, $u_0 \in \bar{S}$ with $\int_0^L \varphi_0(x)dx = \int_0^L u_0(x)dx = 0$ and

$$\varphi_0'''(0) = \varphi_0'''(L) = -\sigma^{-1}K\left(k_5^2 w_0 + k_5\xi_0 - h^*(k_3+k_4)\int_0^L x u_0(x)dx - k_3\mu h^*\left(\varphi_0(L) - \varphi_0(0)\right)\right)$$

be given. We notice that the feedback law (5.6) takes the form (5.5) with

$$\begin{aligned} &B = Kk_5\ ,\ C = Kk_5^2 \\ &\tilde{r}(x) = -Kh^*(k_3+k_4)x \\ &\tilde{p}(x) = -Kk_3\mu h^* \end{aligned} \tag{6.54}$$

For $\phi_n(x) = \sqrt{\frac{2}{L}}\cos\left(n\pi\frac{x}{L}\right)$, $n = 1,2,...$, we get from (6.54) for $n = 1,2,...$:

$$\begin{aligned} \langle \phi_n, \tilde{r} \rangle &= -Kh^*(k_3+k_4)\frac{L^2\left((-1)^n - 1\right)}{n^2\pi^2}\sqrt{\frac{2}{L}} \\ \langle \phi_n', \tilde{p} \rangle &= -Kk_3\mu h^*\sqrt{\frac{2}{L}}\left((-1)^n - 1\right) \end{aligned} \tag{6.55}$$

Using (5.4), (6.54) and (6.55) we get:

$$\begin{aligned}
&\langle \tilde{g}, \tilde{r} \rangle = Kh^*(k_3 + k_4)\frac{L^5}{20} \\
&\frac{12L^4}{\pi^4}\sqrt{\frac{2}{L}}\sum_{n\,odd}\frac{\left|\langle \phi_n, \tilde{r} \rangle\right|}{n^4} = \frac{48L^5}{\pi^6}Kh^*(k_3 + k_4)\sum_{n\,odd}\frac{1}{n^6} \\
&\langle \tilde{g}', \tilde{p} \rangle = Kk_3\mu h^*\frac{L^3}{2} \\
&\frac{12L^4}{\pi^4}\sqrt{\frac{2}{L}}\sum_{n\,odd}\frac{\left|\langle \phi_n', \tilde{p} \rangle\right|}{n^4} = \frac{48L^3}{\pi^4}Kk_3\mu h^*\sum_{n\,odd}\frac{1}{n^4}
\end{aligned} \tag{6.56}$$

Since $\sum_{n\,odd}\frac{1}{n^6} = \frac{\pi^6}{960}$ and $\sum_{n\,odd}\frac{1}{n^4} = \frac{\pi^4}{96}$, it follows from (6.56) that inequalities (5.3) hold. Therefore, Theorem 6 implies that there exist unique functions $\varphi \in C^0\left(\mathbb{R}_+; \bar{S} \cap H^4(0,L)\right) \cap C^1\left(\mathbb{R}_+; \bar{S}\right) \cap C^2\left(\mathbb{R}_+; L^2(0,L)\right)$, $(\xi, w) \in C^1\left(\mathbb{R}_+; \mathbb{R}^2\right)$ with $(\xi(0), w(0)) = (\xi_0, w_0)$, $\varphi[0] = \varphi_0$, $\varphi_t[0] = u_0$ that satisfy (4.1)-(4.5) and (5.6) for all $t \geq 0$.

The rest of the proof exploits the Lyapunov functional

$$\begin{aligned}
\tilde{W}(t) &= \frac{1}{2}\xi^2(t) + \frac{k_5^2}{2}\left(w(t) + k_5^{-1}\xi(t)\right)^2 \\
&+ \frac{\mu}{K\left(h^*\right)^2 L}\left(\frac{1}{2}\|\varphi_t[t]\|_2^2 + \frac{c^2}{2}\|\varphi_x[t]\|_2^2 + \frac{\sigma h^*}{2}\|\varphi_{xx}[t]\|_2^2\right) \\
&+ k_4\left(\frac{1}{2}\|\theta[t]\|_2^2 + \frac{c^2}{2}\|\varphi[t]\|_2^2 + \frac{\sigma h^*}{2}\|\varphi_x[t]\|_2^2\right) \\
&+ k_3\left(\frac{1}{2}\|\theta[t] - \mu\varphi_x[t]\|_2^2 + \frac{c^2 + \bar{\kappa}\mu}{2}\|\varphi[t]\|_2^2 + \frac{\sigma h^*}{2}\|\varphi_x[t]\|_2^2\right)
\end{aligned} \tag{6.57}$$

where $\theta$ is defined by (6.18). Using (6.21), (6.22), (6.23) and definition (6.57) we get for all $t \geq 0$:

$$\begin{aligned}
\frac{d}{dt}\left(\tilde{W}(t)\right) &= -k_5^{-1}\xi^2(t) + k_5\left(w(t) + k_5^{-1}\xi(t)\right)^2 - \frac{\mu^2}{K\left(h^*\right)^2 L}\|\varphi_{tx}[t]\|_2^2 \\
&- \left(\frac{\bar{\kappa}\mu}{K\left(h^*\right)^2 L} + k_4\mu\right)\|\varphi_t[t]\|_2^2 + \frac{\mu}{Kh^*L}f(t)\int_0^L \varphi_{tx}(t,x)dx \\
&- \bar{\kappa}\left(k_3 + k_4\right)\|\theta[t]\|_2^2 - k_3\mu c^2\|\varphi_x[t]\|_2^2 - k_3\mu\sigma h^*\|\varphi_{xx}[t]\|_2^2 \\
&- \left(k_5^2 w(t) + k_5\xi(t) + h^*(k_3 + k_4)\int_0^L \theta(t,x)dx - k_3\mu h^*\int_0^L \varphi_x(t,x)dx\right)f(t)
\end{aligned} \tag{6.58}$$

Using integration by parts and (6.18), (6.20) we get from (5.6) for all $t \geq 0$:

$$f(t) = K\left(k_5^2 w(t) + k_5\xi(t) + h^*(k_3 + k_4)\int_0^L \theta(t,x)dx - k_3\mu h^*\int_0^L \varphi_x(t,x)dx\right) \tag{6.59}$$

Combining (6.58) and (6.59) and using the fact that $\left|\int_0^L \varphi_{tx}(t,x)dx\right| \le \sqrt{L}\,\|\varphi_{tx}[t]\|_2$, we get for all $t \ge 0$ and $\bar{\gamma} \ge 0$:

$$\begin{aligned}
&\frac{d}{dt}\left(\tilde{W}(t)\right) \le -k_5^{-1}\xi^2(t) - \bar{\gamma}\left(w(t) + k_5^{-1}\xi(t)\right)^2 - K^{-1}f^2(t) \\
&+\left(k_5 + \bar{\gamma}\right)\left(w(t) + k_5^{-1}\xi(t)\right)^2 - \frac{\mu^2}{K\left(h^*\right)^2 L}\|\varphi_{tx}[t]\|_2^2 \\
&-\left(\frac{\bar{\kappa}\mu}{K\left(h^*\right)^2 L} + k_4\mu\right)\|\varphi_t[t]\|_2^2 + \frac{\mu}{Kh^*\sqrt{L}}|f(t)|\|\varphi_{tx}[t]\|_2 \\
&-\bar{\kappa}\left(k_3 + k_4\right)\|\theta[t]\|_2^2 - k_3\mu c^2\|\varphi_x[t]\|_2^2 - k_3\mu\sigma h^*\|\varphi_{xx}[t]\|_2^2
\end{aligned} \tag{6.60}$$

Combining (6.59) and (6.60) and using the fact that $\frac{\mu}{Kh^*\sqrt{L}}|f(t)|\|\varphi_{tx}[t]\|_2 \le \frac{\mu^2}{2K\left(h^*\right)^2 L}\|\varphi_{tx}[t]\|_2^2 + \frac{1}{2K}f^2(t)$, we get for all $t \ge 0$ and $\bar{\gamma} \ge 0$:

$$\begin{aligned}
&\frac{d}{dt}\left(\tilde{W}(t)\right) \le -k_5^{-1}\xi^2(t) - \bar{\gamma}\left(w(t) + k_5^{-1}\xi(t)\right)^2 \\
&+k_5^{-3}\left(1 + k_5^{-1}\bar{\gamma}\right)\left(K^{-1}f(t) - h^*(k_3 + k_4)\int_0^L \theta(t,x)dx + k_3\mu h^*\int_0^L \varphi_x(t,x)dx\right)^2 \\
&-\frac{\mu^2}{2K\left(h^*\right)^2 L}\|\varphi_{tx}[t]\|_2^2 - \left(\frac{\bar{\kappa}\mu}{K\left(h^*\right)^2 L} + k_4\mu\right)\|\varphi_t[t]\|_2^2 - \frac{1}{2K}f^2(t) \\
&-\bar{\kappa}\left(k_3 + k_4\right)\|\theta[t]\|_2^2 - k_3\mu c^2\|\varphi_x[t]\|_2^2 - k_3\mu\sigma h^*\|\varphi_{xx}[t]\|_2^2
\end{aligned} \tag{6.61}$$

Using the fact that

$$\begin{aligned}
&\left(K^{-1}f(t) - h^*(k_3 + k_4)\int_0^L \theta(t,x)dx + k_3\mu h^*\int_0^L \varphi_x(t,x)dx\right)^2 \\
&\le 2K^{-2}f^2(t) + 4k_3^2\mu^2(h^*)^2\left(\int_0^L \varphi_x(t,x)dx\right)^2 + 4(h^*)^2(k_3 + k_4)^2\left(\int_0^L \theta(t,x)dx\right)^2
\end{aligned}$$

we obtain from (6.61) for all $t \ge 0$ and $\bar{\gamma} \ge 0$:

$$\frac{d}{dt}\left(\tilde{W}(t)\right) \le -k_5^{-1}\xi^2(t) - \bar{\gamma}\left(w(t) + k_5^{-1}\xi(t)\right)^2$$
$$+4k_5^{-3}\left(1+k_5^{-1}\bar{\gamma}\right)k_3^2\mu^2(h^*)^2\left(\int_0^L \varphi_x(t,x)dx\right)^2 + 4k_5^{-3}\left(1+k_5^{-1}\bar{\gamma}\right)(h^*)^2(k_3+k_4)^2\left(\int_0^L \theta(t,x)dx\right)^2$$
$$-\frac{\mu^2}{2K\left(h^*\right)^2 L}\|\varphi_{tx}[t]\|_2^2 - k_4\mu\|\varphi_t[t]\|_2^2 - 2K^{-2}\left(\frac{K}{4} - k_5^{-3}\left(1+k_5^{-1}\bar{\gamma}\right)\right)f^2(t) \tag{6.62}$$
$$-\bar{\kappa}\left(k_3+k_4\right)\|\theta[t]\|_2^2 - k_3\mu c^2\|\varphi_x[t]\|_2^2 - k_3\mu\sigma h^*\|\varphi_{xx}[t]\|_2^2$$

Using the facts that $\left(\int_0^L \varphi_x(t,x)dx\right)^2 \le L\|\varphi_x[t]\|_2^2$ and $\left(\int_0^L \theta(t,x)dx\right)^2 \le L\|\theta[t]\|_2^2$, we obtain from (6.62) for all $t \ge 0$ and $\bar{\gamma} \ge 0$:

$$\frac{d}{dt}\left(\tilde{W}(t)\right) \le -k_5^{-1}\xi^2(t) - \bar{\gamma}\left(w(t) + k_5^{-1}\xi(t)\right)^2$$
$$+4k_5^{-3}\left(1+k_5^{-1}\bar{\gamma}\right)(h^*)^2(k_3+k_4)^2 L\|\theta[t]\|_2^2 - k_3\mu\sigma h^*\|\varphi_{xx}[t]\|_2^2$$
$$-\frac{\mu^2}{2K\left(h^*\right)^2 L}\|\varphi_{tx}[t]\|_2^2 - k_4\mu\|\varphi_t[t]\|_2^2 - 2K^{-2}\left(\frac{K}{4} - k_5^{-3}\left(1+k_5^{-1}\bar{\gamma}\right)\right)f^2(t) \tag{6.63}$$
$$-4k_3^2\mu^2(h^*)^2 L\left(\frac{c^2}{4k_3\mu(h^*)^2 L} - k_5^{-3}\left(1+k_5^{-1}\bar{\gamma}\right)\right)\|\varphi_x[t]\|_2^2$$

Using (6.27) and (6.29), we obtain from (6.63) for all $t \ge 0$ and $\bar{\gamma} \ge 0$:

$$\frac{d}{dt}\left(\tilde{W}(t)\right) \le -k_5^{-1}\xi^2(t) - \bar{\gamma}\left(w(t) + k_5^{-1}\xi(t)\right)^2$$
$$-k_3\mu\sigma h^*\|\varphi_{xx}[t]\|_2^2 - 2K^{-2}\left(\frac{K}{4} - k_5^{-3}\left(1+k_5^{-1}\bar{\gamma}\right)\right)f^2(t)$$
$$-4(h^*)^2(k_3+k_4)^2\frac{L^3}{\pi^2}\left(\frac{\mu\pi^2\left(2Kk_4(h^*)^2L^3 + \mu\pi^2\right)}{8K(h^*)^4(k_3+k_4)^2L^6} - k_5^{-3}\left(1+k_5^{-1}\bar{\gamma}\right)\right)\|\varphi_t[t]\|_2^2 \tag{6.64}$$
$$-4k_3^2\mu^2(h^*)^2 L\left(\frac{c^2}{4k_3\mu(h^*)^2 L} - k_5^{-3}\left(1+k_5^{-1}\bar{\gamma}\right)\right)\|\varphi_x[t]\|_2^2$$

Inequality (5.7) implies that there exists $\bar{\gamma} > 0$ sufficiently small such that $K > 4k_5^{-3}\left(1+k_5^{-1}\bar{\gamma}\right)$, $\frac{\mu\pi^2\left(2Kk_4(h^*)^2L^3 + \mu\pi^2\right)}{8K(h^*)^4(k_3+k_4)^2L^6} > k_5^{-3}\left(1+k_5^{-1}\bar{\gamma}\right)$ and $\frac{c^2}{4k_3\mu(h^*)^2 L} > k_5^{-3}\left(1+k_5^{-1}\bar{\gamma}\right)$. Consequently, there exists a constant $c_1 > 0$ (independent of the solution and independent of $t \ge 0$) such that the following inequality holds for all $t \ge 0$:

$$\frac{d}{dt}\left(\tilde{W}(t)\right) \le -c_1\left(\xi^2(t) + \left(w(t) + k_5^{-1}\xi(t)\right)^2 + \|\varphi_x[t]\|_2^2 + \|\varphi_{xx}[t]\|_2^2 + \|\varphi_t[t]\|_2^2\right) \tag{6.65}$$

Using definition (6.57) and inequalities (6.27), (6.28) we conclude that there exist a constant $c_2 > 0$ (independent of $t \geq 0$ and the solution) such that the following inequality holds for all $t \geq 0$:

$$\tilde{W}(t) \leq K_1 \left( \xi^2(t) + \left( w(t) + k_5^{-1} \xi(t) \right)^2 + \|\varphi_x[t]\|_2^2 + \|\varphi_{xx}[t]\|_2^2 + \|\varphi_t[t]\|_2^2 \right) \tag{6.66}$$

Combining (6.65) and (6.66) we obtain the differential inequality $\frac{d}{dt}\left(\tilde{W}(t)\right) \leq -c_2^{-1} c_1 \tilde{W}(t)$ for all $t \geq 0$ which directly gives the following estimate for all $t \geq 0$:

$$\tilde{W}(t) \leq \exp\left(-K_1^{-1} R t\right) \tilde{W}(0) \tag{6.67}$$

The rest of the proof follows from estimate (6.67) and the fact that there exist constants $c_4 > c_3 > 0$ (independent of $t \geq 0$ and the solution $\varphi, \xi, w$) such that the following inequalities hold for all $t \geq 0$:

$$c_3 \left( P^2(t) + \xi^2(t) + w^2(t) \right) \leq \tilde{W}(t) \leq c_4 \left( P^2(t) + \xi^2(t) + w^2(t) \right) \tag{6.68}$$

The proof is complete. $\lhd$

## 7. Open Problems

The present paper reviewed the available results for the feedback stabilization problem of the viscous tank-liquid system. We have also given novel results for the linearization of the tank-liquid system and we have showed that the same family of feedback laws that works for the nonlinear case can also achieve exponential stabilization for the linearization of the tank-liquid system, no matter what the value of the surface tension coefficient is.

However, many new results are still needed for a complete study of the spill-free, slosh-free problem of tank-liquid transfer. The following list presents a number of issues that remain unresolved.

<u>1) Extension to two dimensions</u>
Real tanks are not one-dimensional; the actual equations for the study of a real tank-liquid system involve an additional spatial dimension. The extension of the theoretical results to two spatial dimensions pose new challenges that will demand different approaches that have no analogues in the 1-D case (e.g., the appearance of vorticity). The extension may also involve different boundary conditions when surface tension is present (see [38, 46]).

<u>2) Existence/Uniqueness</u>
As noted above, the nonlinear results of Section 3 are not accompanied by appropriate existence/uniqueness results for the solutions of the closed-loop system. The effort for the development of new existence/uniqueness results will necessarily involve weaker solution notions than the ones that were used in Section 3. Therefore, this effort will also involve the extension of the existing stabilization results to weaker solution notions and will require new mathematical tools (see the beginning of such an effort in [37]).

3) Addressing the nonlinear case where both surface tension and friction are present

It should be emphasized that the nonlinear case where both surface tension and friction has not been addressed so far. The results of Section 5 for the linearization of the tank-liquid system give hopes that the nonlinear problem is solvable. Such a research effort will also require the derivation of additional stability estimates that provide bounds for the sup-norm of the velocity in the presence of surface tension.

4) The need for numerical methods

None of the numerical schemes proposed in the literature can be applied to the case where the tank-liquid system is under the effect of a feedback controller. This is true not only for the nonlinear tank-liquid PDE-ODE system but also for the linearized tank-liquid system. Therefore, the development of novel numerical schemes is needed for the accurate and reliable simulation of the closed-loop system. Such a research effort has already started in [34, 37] where promising numerical schemes that preserve the stability estimates have been developed. It should be also noted that the need of fast and accurate numerical schemes is very important for the control practitioner who will implement the feedback controllers. The practical implementation of the proposed controllers involves numerous numerical experiments for the tuning of the controller parameters and the evaluation of their performance.

5) Adaptive controllers

The actual value of some of the constants $\mu, g, \sigma, h^*, L$ may be unknown. Therefore, adaptive controllers need to be developed in order to address such a case.

6) Output feedback stabilization

The need for output feedback stabilizers comes from the difficulty to obtain measurements for all the quantities that are used in the proposed feedback laws (e.g., the liquid momentum or the tank velocity). The study of observer-based output feedback laws started recently in [34] but needs to be extended to the case where friction and/or surface tension is present.

7) Rapid Stabilization (at least in the linearized case)

It would be interesting to construct feedback laws that can achieve rapid stabilization of the tank-liquid system, i.e., stabilization with an (arbitrarily) assignable, exponential convergence rate (see [10] for the rapid stabilization of beam equations and see [54] for stabilization of linearized viscous flows). Notice that the existing stabilization results do not guarantee an (arbitrarily) assignable, exponential convergence rate. If the rapid stabilization problem is not solvable, it would also be interesting to have a counterexample that shows that rapid stabilization is not possible.